\documentclass[3p,11pt]{elsarticle}
\usepackage{lineno,hyperref}
\modulolinenumbers[1]
\usepackage{geometry}
\usepackage{graphicx}
\usepackage[usenames,dvipsnames]{color}
\usepackage{color}
\usepackage[english]{babel}
\usepackage{epstopdf}
\usepackage{latexsym}
\usepackage{hyperref}
\usepackage{url}
\usepackage{amstext}
\usepackage{amssymb}
\usepackage{amsmath}
\usepackage{pifont}
\usepackage{url}
\hypersetup{
    colorlinks=true,
    linkcolor=blue,
    filecolor=magenta,
    urlcolor=blue,
}
\graphicspath{ {Images/} }
\usepackage{times}
\usepackage{hyperref}
\usepackage{url}
\usepackage{amstext}
\usepackage{amssymb}
\usepackage{wrapfig}
\usepackage{float}
\usepackage{tikz}
\usepackage{pgfplots}
\usepackage{caption}
\usepackage{subcaption}

\begin{document}

\begin{frontmatter}

\title{Stationary and non-stationary pattern formation over fragmented habitat}

\author[inst1]{Malay Banerjee\corref{cor1}}
\ead{malayb@iitk.ac.in}
\author[inst1,inst2]{Swadesh Pal}
\ead{palkkeswadesh@gmail.com}
\author[inst1]{Pranali Roy Chowdhury}
\ead{pranali@iitk.ac.in}

\cortext[cor1]{Corresponding author}
\address[inst1]{Department of Mathematics and Statistics,
IIT Kanpur, Kanpur, India}
\address[inst2]{MS2Discovery Interdisciplinary Research Institute, Wilfrid Laurier University, Waterloo, Canada}

\begin{abstract}
Spatio-temporal pattern formation over the square and rectangular domain has received significant attention from researchers. A wide range of stationary and non-stationary patterns produced by two interacting populations is abundant in the literature. Fragmented habitats are widespread in reality due to the irregularity of the landscape. This work considers a prey-predator model capable of producing a wide range of stationary and time-varying patterns over a complex habitat. The complex habitat is assumed to have consisted of two rectangular patches connected through a corridor. Our main aim is to explain how the shape and size of the fragmented habitat regulate the spatio-temporal pattern formation at the initial time. The analytical conditions are derived to ensure the existence of a stationary pattern and illustrate the role of most unstable eigenmodes to determine the number of patches for the stationary pattern. Exhaustive numerical simulations help to explain the spatial domain's size and shape on the transient patterns and the duration of transient states.
\end{abstract}


\end{frontmatter}

\section{Introduction}
The complex interaction among ecological species is ubiquitous and is intrinsically fascinating. This has been a key focus of research over many decades, where the applied mathematicians worked hand in hand with the ecologists \cite{Kar90,Ros63,Tur03}. The spatial distribution of the species population studied over a specific region with the help of reaction-diffusion equations \cite{Cantrell03}. In this context, the seminal work of Turing \cite{Tur52} and later by Segel and Jackson \cite{Segel72} gave a new insight into spatio-temporal pattern formation to understand the spatial heterogeneity of species in the ecosystem effectively. The mathematical models consisting of reaction-diffusion equations provide a valuable framework to investigate self-organized pattern formation. But, for the tractability of mathematical models, most studies considered the spatial distribution of species in uniform domains. However, this is far from reality. Consequently, over the last few years, considerable attention was given to understand the spatial spread of species and thus the pattern formation in much more complex domains \cite{Alharbi15,Alharbi18,Alharbi19}.

The mobility of the species from one region to another depends on various environmental factors, including climate condition, habitat fragmentation, species competition, etc \cite{Malvido08,Ewers01}. It has a massive impact on the ecosystem. One of the major concerns is the fragmentation of the ecological habitats, which was considered an invasive threat to biodiversity. Habitat fragmentation can define as a landscape-scale process in which the continuous habitat is reduced into smaller habitat remnants \cite{Malvido08}. The size, shape, edge of the habitat fragments, and habitat isolation are some significant factors having huge implications on the species interaction and species survival \cite{Alharbi15, Ewers01}. Researchers have studied these implications over different temporal scales, with short term effects including change in population sizes or edge effects of the fragments and long term effects dealing with change in the pattern of gene flow or extinction of species \cite{Alharbi15,Browne15}.

Although habitat fragmentation was an active area of research for a long time in ecology, very few works are available in mathematical modelling. In \cite{Alharbi15}, the authors explore how the size, shape, and edge of the fragment habitats influence the population dynamics of the species. They analyzed it in three different domains, namely, square, cross-shaped, and H-shaped. Considering the Allee effect, they have estimated the critical patch size for the species' existence. They showed that if the critical patch size of the cross-shaped domain is larger than the corresponding square-shaped domain, then the likelihood of population survival is higher in the simple domain. Also, with an increase in strength of the Allee effect, the critical size increase. The inclusion of biological corridors between population patches seemed to have beneficial effects on the movement of the species and thus on the persistence of the species \cite{Ewers01}. In \cite{Alharbi18,Alharbi19}, researchers have studied the effects of the corridor on the invasion of alien species, more specifically, how the length and width of the corridor alter the speed of invasion in an H-shaped domain. Initially, they have introduced an alien species at a particular spatial point in one of the habitats and showed that the propagation of the population front to the other habitat is considerably slow for a smaller width of the corridor. The corridor can also modify the pattern of spread. In the case of patchy invasion in one of the habitats, the movement of the patches is random. Suppose the corridor's length and width are greater than a critical value. In that case, to accommodate the random patches to penetrate the other habitat, the patchy invasion occurs in the entire domain. Sometimes, the corridor can regularize the spatio-temporal pattern formation. We mean to say the irregular spatio-temporal pattern in one habitat changes to regular circular fronts in the other habitat. Later, the authors considered a small patch in the corridor with a better environmental condition called a stepping stone and showed its importance in the dispersal and persistence of the species. They investigated the effect of shape, size, location, and quality of the stepping stone on the spread of the alien species in fragmented habitats.

We observe the spatio-temporal interaction between prey and predator's influence on the spatial patterns in nature. In this context, the mathematical framework of diffusion-driven instability or Turing bifurcation helps to study various pattern formations. Under the conditions of Turing bifurcation, the homogeneous steady state is stable in the absence of diffusion, whereas unstable in the presence of diffusion. The patterns arising because of this instability are known as Turing patterns. Turing patterns classify as spots (hot or cold), stripes, labyrinthine, and a mixture of stripes and spots. The prey-predator models with a prey-dependent functional response and linear death rate of the predator cannot support Turing patterns. However, considering predator interference in the functional response, the nonlinear death rate of the predator or intra-specific competition among predators can lead to Turing patterns.

The main objective of this paper is to study the effect of fragmented habitat on spatio-temporal pattern formation. Particularly, we consider a spatio-temporal prey-predator model with additive weak Allee effect in prey growth in a U-shaped domain. In this work, we focus on different stationary and dynamic pattern formation in a U-shaped domain, where two large habitats are connected by a narrow corridor. This corridor facilitates the movement of the individuals from one habitat to the other, and thus the pattern formation in both the habitats. The paper is organized as follows. We present the spatio-temporal model in Section 2 and discuss the analytical conditions of local stability, Hopf bifurcation, and Turing bifurcation. In Section 3, we compare the different pattern formation scenarios in the square and U-shaped domains with the help of exhaustive numerical simulations. We show how the size of the habitats and the connecting channel between the habitats affects the spatio-temporal pattern formation for different initial conditions. We further explore the transients dynamics observed in the fragmented habitat and how the connecting channel's width affects the transient time and amplitude of the oscillations. Finally, we conclude our work in Section 4.

\section{Mathematical model}

In this work, we consider the following two-dimensional spatio-temporal model of prey-predator interaction with Bazykin type reaction kinetics and Allee effect in prey growth \cite{ WWIJB,MBWWMBE, KMMBEC, KMMBMBE}:
\begin{subequations}{\label{model1}}
\begin{align}
\frac{\partial u}{\partial t} & = d_{1}\bigg{(}\frac{\partial^{2} u}{\partial x^{2}}+\frac{\partial^{2} u}{\partial y^{2}}\bigg{)} + u\bigg{(} r-fu-\frac{m}{b+u}-\frac{cv}{u+a}\bigg{)}, \\
\frac{\partial v}{\partial t} & = d_{2}\bigg{(}\frac{\partial^{2} v}{\partial x^{2}}+\frac{\partial^{2} v}{\partial y^{2}}\bigg{)} + sv\bigg{(} \frac{cu}{u+a}-q-pv\bigg{)},
\end{align}
\end{subequations}
subjected to non-negative initial conditions and zero-flux boundary conditions. Here, we consider mainly the U-shaped domain as an example of fragmented habitat to understand the spatio-temporal pattern formation over a non-conventional spatial domain. All the parameters involved with the model are positive constants. The parameters $r$ and $f$ denote  the intrinsic growth rate and intra-specific competition strength of the prey population. The additive Allee effect is modelled by the term $bu/(m+u)$, where $b$ measures the strength of the Allee effect and $m$ is the prey population size at which the fitness is half of the maximum value. Consumption of prey by their specialist predator follows saturating functional response, $c$ is the prey capturing rate by the predator and $a$ measures the half-saturation level. The predator's natural death rate and density-dependent death rate are denoted by $q$ and $p$, respectively. The multiplicative parameter $s$ is the feed concentration for the specialist predator. In a crude sense, the growth of the predator population follows a logistic type growth law $(dv/dt=\alpha v -\beta v^2)$ where the growth rate is prey dependent $\alpha\equiv csu/(u+a)-sq$ and $\beta=sp$.

The additional term $bu/(m+u)$ is introduced in the prey growth rate through the addition of a negative feedback term to model the Allee effect, and it is known as the additive Allee effect \cite{GonzalezSIAM,GonzalezNARWA}. The phrase `additive Allee effect' is nothing to do with the ecology; instead, some parametric restrictions determine the ecological aspect. Without going into the mathematical details, as the results are well known, we can mention that the Allee effect is strong when the parameters satisfy the restriction
\begin{eqnarray}\label{StrongAllee}
b^2rf<br<m<\frac{r^2(1-bf)^2+4br^2f}{4rf}
\end{eqnarray}
and the Allee effect is weak under the parametric restriction
\begin{eqnarray}\label{WeakAllee}
b^2rf<m<br.
\end{eqnarray}
The temporal model corresponding to the spatio-temporal model (\ref{model1}) is given by
\begin{subequations}{\label{model2}}
\begin{align}
\frac{du}{dt} & = u\bigg{(} r-fu-\frac{m}{b+u}-\frac{cv}{u+a}\bigg{)}\,\equiv\,ug(u)-h(u)v, \\
\frac{dv}{dt} & = sv\bigg{(} \frac{cu}{u+a}-q-pv\bigg{)}\,\equiv\,s\left(h(u)-m(v)\right)v,
\end{align}
\end{subequations}
with non-negative initial conditions and the functions $g(.)$, $h(.)$ and $m(.)$ are given by
$$g(u)\,=\,r-fu-\frac{m}{b+u},\,\,h(u)\,=\,\frac{cu}{u+a},\,\,m(v)\,=\,q+pv.$$
These three functions represent the per capita growth rate of the prey in the absence of a predator, prey-dependent functional response and density-dependent per capita predator's death rate, respectively. The model under consideration is a Gause type prey-predator model with specialist predator and additive Allee effect in prey growth.

As per the primary goal of this work, we consider the case of a unique coexisting equilibrium point, which is the unique point of intersection of non-trivial nullclines $g(u)-vh(u)/u=0$ and $h(u)=m(v)$ in the interior of the first quadrant. Let $E_*(u_*,v_*)$ be the non-trivial equilibrium point, then $v_*=u_*g(u_*)/h(u_*)$, $u_*=h^{-1}(m(v_*))$. In other words, $u_*$ is a positive root of the quadratic equation
$$A_0u^4+A_1u^3+A_2u^2+A_3u+A_4=0,$$
where $A_0=rpf$, $A_1=rp((2a+b)f-1)$, $A_2=a(a+2b)rpf+mp+c(c-q)-(2a+b)rp$, $A_3=a^2brpf+2amp+bc (c-q)-a(a + 2b)rp- acq$, $A_4=a^2p(m-br)-abcq$, and $v_*$ is given by
$$v_*=\frac{1}{p}\left(\frac{cu_*}{u_*+a}-q\right).$$
For the feasibility of $v_*$, we require the restrictions $u_*>ap/(c-q)$ and $c>a$.

\subsection{Local stability and Hopf-bifurcation}

We use Routh-Hurwitz criteria to determine the local asymptotic stability of coexisting equilibrium point \cite{Perko}. The Jacobian matrix of the system (\ref{model2}) evaluated at $E_*$ is given by
\begin{eqnarray}\label{Jacobian}
J_*\,=\,\left[\begin{array}{cc}
g(u_*)+u_*g'(u_*)-h'(u_*)v_* & -h(u_*) \\
sh'(u_*)v_* & -spv_* \\
\end{array}\right]\,\equiv\,\left[\begin{array}{cc}
\alpha_{11} & \alpha_{12} \\
\alpha_{21} & \alpha_{22} \\
\end{array}\right].
\end{eqnarray}
According to the Routh-Hurwitz criteria \cite{HirschSmale}, $E_*$ is locally asymptotically stable if the following two conditions are satisfied:
\begin{eqnarray}\label{stab}
\textrm{tr}(J_*)\,=\,\alpha_{11}+\alpha_{22}\,<\,0,\,\,\textrm{det}(J_*)\,=\,\alpha_{11}\alpha_{22}-\alpha_{12}\alpha_{21}\,>\,0.
\end{eqnarray}
Components of $E_*$ are positive, parameters are positive and $h'(u)>0$ for all $u>0$. Except $\alpha_{11}$, signs of other entries of $J_*$ are fixed and hence the sufficient condition for local asymptotic stability of $E_*$ is given by $\alpha_{11}<0$.

The coexisting steady-state $E_*$ loses stability through Hopf bifurcation whenever $\alpha_{11}>0$. The Hopf bifurcation condition for the coexistence equilibrium $E_*$ is given by
\begin{eqnarray}\label{hopf}
\textrm{tr}(J_*)\,=\,0,\,\,\textrm{det}(J_*)\,>\,0,\,\,\frac{d}{ds}\left[\textrm{tr}(J_*)\right]\,\neq\,0.
\end{eqnarray}
It is interesting to mention that the components of $E_*$ can not be found explicitly, however, the model formulation allow us to obtain the Hopf-bifurcation threshold explicitly in terms of $s$ as:
\begin{eqnarray}\label{hopf1}
s_H=\frac{g(u_*)+u_*g'(u_*)-h'(u_*)v_*}{pv_*}.
\end{eqnarray}
Note that the Hopf bifurcation threshold is feasible when $\alpha_{11}>0$ and $\alpha_{11}\alpha_{22}>\alpha_{12}\alpha_{21}$. One can find a limit cycle in the vicinity of Hopf bifurcation threshold, and the stability of Hopf bifurcating limit cycle depends on the sign of the first Lyapunov number \cite{Perko,Kuznetsov04}. To obtain the expression for the first Lyapunov number, we first expand the Taylor series expansion of the system (\ref{model2}) around $E_*$, $s=s_H$ as the following
\begin{subequations}\label{taylor}
\begin{eqnarray}
\frac{du_1}{dt}&=&\beta_{10}u_1+\beta_{01}v_1+\beta_{20}u_1^2+\beta_{11}u_1v_1+\beta_{02}v_2^2+\beta_{30}u_1^3+\beta_{21}u_1^2v_1+\beta_{12}u_1v_1^2+\beta_{03}v_2^3+\cdots,\\
\frac{dv_1}{dt}&=&\gamma_{10}u_1+\gamma_{01}v_1+\gamma_{20}u_1^2+\gamma_{11}u_1v_1+\gamma_{02}v_2^2+\gamma_{30}u_1^3+\gamma_{21}u_1^2v_1+\gamma_{12}u_1v_1^2+\gamma_{03}v_2^3+\cdots.
\end{eqnarray}
\end{subequations}
The expression for the first Lyapunov number is given by
\begin{eqnarray}\label{hopf2}
\sigma&=&-\frac{3\pi}{2\beta_{01}\Delta^{3/2}}\left[
\beta_{10}\gamma_{10}(\beta_{11}^2+\beta_{11}\gamma_{02}+\beta_{02}\gamma_{11})+\beta_{10}\beta_{01}(\gamma_{11}^2+\beta_{20}\gamma_{11}+\beta_{11}\gamma_{02})\right.\nonumber\\
& & +\gamma_{10}^2(\beta_{11}\beta_{02}+2\beta_{02}\gamma_{02})-2\beta_{10}\gamma_{10}(\gamma_{02}^2-\beta_{02}\beta_{20})-2\beta_{10}\beta_{01}(\beta_{20}^2-\gamma_{20}\gamma_{02})\nonumber\\
&& -\beta_{01}^2(2\beta_{20}\gamma_{20}+\gamma_{11}\gamma_{20})+(\beta_{01}\gamma_{10}-2\beta_{10}^2)(\gamma_{11}\gamma_{02}-\beta_{11}\beta_{20})\nonumber\\
&&\left.-(\beta_{10}^2+\beta_{01}\gamma_{10})(3(\gamma_{10}\gamma_{03}-\beta_{01}\beta_{30})+2\beta_{10}(\beta_{21}+\gamma_{12})+(\gamma_{10}\beta_{12}-\beta_{01}\gamma_{21}))\right],
\end{eqnarray}
where $\Delta=\beta_{10}\gamma_{01}-\beta_{01}\gamma_{10}$. The Hopf bifurcation is super-critial if $\sigma<0$ and is sub-critical for $\sigma>0$. We verify the stability of Hopf bifurcating limit cycle with the help of numerical example.

\subsection{Turing bifurcation}

Here we derive the Turing instability condition for the system (\ref{model1}) around the homogeneous steady-state $E_*$, where $u(t,x,y)=u_*$, $v(t,x,y)=v_*$ is the homogeneous steady-state for (\ref{model1}). The model (\ref{model1}) is subject to the non-negative initial condition
\begin{equation}\label{BC}
    u(0,x,y)\,=\,u_0(x,y)\,\geq\,0,\,\,v(0,x,y)\,=\,v_0(x,y)\,\geq\,0,\,\,0\,\leq\,x,y\,\leq\,L,
\end{equation}
and zero-flux boundary condition along the boundary of the square domain.

After linearizing the system (\ref{model1}) around the homogeneous steady-state $E*$, we obtain the following system of linear PDEs in terms of the perturbation variables $u_1(t,x,y)$ and $v_1(t,x,y)$:
\begin{equation}\label{LinSTeq}
\begin{aligned}
      \dfrac{\partial u_1}{\partial t} &= d_{1}\bigg{(}\dfrac{\partial ^2u_1}{\partial x^2}+\dfrac{\partial ^2u_1}{\partial y^2}\bigg{)} + \alpha_{11}u_1+\alpha_{12}v_1,\\
      \dfrac{\partial v_1}{\partial t} &= d_{2}\bigg{(}\dfrac{\partial ^2v_1}{\partial x^2}+\dfrac{\partial ^2v_1}{\partial y^2}\bigg{)} + \alpha_{21}u_1+\alpha_{22}v_1.
    \end{aligned}
\end{equation}
We assume that the solution of the system of linear equations (\ref{LinSTeq}) satisfying the no-flux boundary conditions and the solution is in the form
\begin{equation}\label{eq:hetero_pertur}
    \begin{aligned}
    u_1(t,x,y) \, =\, \epsilon_1 e^{\lambda t}\cos(k_1x)\cos(k_2y),\,\,
    v_1(t,x,y) \, =\, \epsilon_2 e^{\lambda t}\cos(k_1x)\cos(k_2y),
  \end{aligned}
\end{equation}
where $\epsilon_1,\epsilon_2\ll 1$ are two arbitrary constants and $k=\sqrt{k_1^2+k_2^2}$ is wave number with $k_1=n_1\pi/L$, $k_2=n_2\pi/L$, $n_1$ and $n_2$ are two natural numbers. For a non-trivial solution (\ref{eq:hetero_pertur}) of the linear system of equations (\ref{LinSTeq}), $\lambda$ satisfies the following characteristics equation
\begin{equation}\label{cheqn}
|M(k)-\lambda I_2|\,=\,0,
\end{equation}
where
\begin{equation}
M(k) = \begin{pmatrix}
    \alpha_{11}-d_1k^2 & \alpha_{12} \\
    \alpha_{21} & \alpha_{22}-d_2k^2
    \end{pmatrix},
\end{equation}
$\alpha_{ij}$'s $(i,j=1,2)$ are given in (\ref{Jacobian}).
After expanding the characteristic equation (\ref{cheqn}), we obtain
\begin{equation}\label{cheqn1}
    \lambda^2-((\alpha_{11}+\alpha_{22})-(d_1+d_2)k^2)\lambda + h(k^2)=0,
\end{equation}
where
\begin{equation}\label{eq:h_k}
    h(k^2) = d_1d_2k^4-(d_2\alpha_{11}+d_1\alpha_{22})k^2+\alpha_{11}\alpha_{22}-\alpha_{12}\alpha_{21}.
\end{equation}

The homogeneous steady-state $E_*$ is stable under small amplitude heterogeneous perturbation if both the eigenvalues of the characteristic equation (\ref{cheqn1}) has a negative real part. Turing instability sets in due to the destabilization of homogeneous steady-state $E_*$, when one root of the characteristic equation (\ref{cheqn1}) passes through zero at some critical wave number. At the same time, the other root of the characteristic equation remains negative. The critical wave number $k_{c}$ can obtain by solving the equation $\frac{dh(k^2)}{dk^2}=0$ for $k^2$, which is given by
\begin{equation}
    k^2_{c} = \dfrac{1}{2d_1d_2}(d_2\alpha_{11}+d_1\alpha_{22}).
\end{equation}
As the critical wave number $k_c$ is a real number, the dimensionless diffusion coefficients $d_1,d_2>0$, feasible existence of $k_c$ demands the satisfaction of the implicit parametric restriction $d_2\alpha_{11}+d_1\alpha_{22}>0$. Turing instability occurs when a stable homogeneous steady-state (which is stable under small amplitude homogeneous perturbation) becomes unstable under small amplitude heterogeneous perturbation. We first assume that the homogeneous steady-state $E_*$ is locally asymptotically stable, i.e., the conditions (\ref{stab}) are satisfied. Hence, for positive $d_1$ and $d_2$, $\alpha_{11}+\alpha_{22}<0$ and $d_2\alpha_{11}+d_1\alpha_{22}>0$ satisfied simultaneously when $\alpha_{11}\alpha_{22}<0$. The equation of Turing bifurcation boundary can be obtained by eliminating $k^2$ from $h(k^2)=0$ with the help of $dh(k^2)/dk^2=0$ as follows,
\begin{equation}\label{Turcurve}
   d_2\alpha_{11}+d_1\alpha_{22} \,=\,2\sqrt{d_1d_2}\sqrt{\alpha_{11}\alpha_{22}-\alpha_{12}\alpha_{21}}.
\end{equation}
The Turing bifurcation curve is feasible under the implicit parametric restrictions $\alpha_{11}\alpha_{22}<0$ and $\alpha_{11}\alpha_{22}>\alpha_{12}\alpha_{21}$. These two conditions can be written in terms of $M(k)$ as follows
\begin{equation}\label{stabcond}
    \textrm{tr}(M(0))\,<\,0,\,\,\textrm{det}(M(0))\,>\,0.
\end{equation}

The wavenumber $k$ is given by $k=\sqrt{\dfrac{n_1^2\pi^2}{L^2}+\dfrac{n_2^2\pi^2}{L^2}}$, where $L$ is the length of the square domain, $n_1$ and $n_2$ are two natural numbers. We use the notation $\omega=\pi/L$. After substituting $k^2=(n_1^2+n_2^2)\omega^2$ in $h(k^2)=0$ and then solving for $d_{2}$, we find the explicit expression for Turing bifurcation boundary corresponding to the admissible unstable eigenmodes $(n_1,n_2)$ as
\begin{equation}
    \begin{aligned}
      d_{2T}(n_1,n_2) \, = \, \dfrac{(n_1^2+n_2^2)\omega^2d_1\alpha_{22}-\alpha_{11}\alpha_{22}+\alpha_{12}\alpha_{21}}{d_1(n_1^2+n_2^2)\omega^2((n_1^2+n_2^2)\omega^2-a_{11})}.
    \end{aligned}
\end{equation}

For a given set of parameter values, suppose the homogeneous coexisting steady-state $E_*$ is stable under temporal perturbation (i.e. the condition (\ref{stabcond}) holds), then $E_*$ is stable under heterogeneous perturbation for $d_2<d_{2T}$ and is unstable for $d_2>d_{2T}$. Depending upon the chosen parameter values involved with the reaction kinetics, one can find more than one set of values for $(n_1,n_2)$ such that $d_{2T}(n_1,n_2)$ is feasible. Consequently, for a fixed set of parameter values, we find multiple feasible thresholds $d_{2T}(n_1,n_2)$, which indicates the existence of multiple heterogeneous stationary solutions. However, it is not easy to determine the number of heterogeneous stationary solutions. We obtain the pattern through numerical simulation is determined by the feasible value of $d_2>d_{2T}(n_1,n_2)$ and the admissible value of $(n_1,n_2)$  corresponding to the maximum real part of $\lambda(k)$ along with the size of the domain and the choice of the initial condition. We will explain this idea in detail with the help of a numerical example in the following section.

The critical Turing bifurcation threshold, in terms of $d_2$, is given by a positive solution of the equation (\ref{Turcurve}) when other parameter values are given, and let us denote this solution as $d_{2T}$. Now if we choose a value of $d_2>d_{2T}$, keeping other parameters fixed, then we can find a range of values of $k$, say $(k_1,k_2)$, such that $h(k^2)<0$ for $k_1<k<k_2$, where $k_1$ and $k_2$ are given by
\begin{subequations}\label{k12}
\begin{eqnarray}
k_1^2&=&\frac{ d_2\alpha_{11}+d_1\alpha_{22}-\sqrt{( d_2\alpha_{11}+d_1\alpha_{22})^2-4d_1d_2(\alpha_{11}\alpha_{22}-\alpha_{12}\alpha_{21})}}{2d_1d_2},\\
k_2^2&=&\frac{ d_2\alpha_{11}+d_1\alpha_{22}+\sqrt{( d_2\alpha_{11}+d_1\alpha_{22})^2-4d_1d_2(\alpha_{11}\alpha_{22}-\alpha_{12}\alpha_{21})}}{2d_1d_2}.
\end{eqnarray}
\end{subequations}
For $d_2>d_{2T}$, one can verify that the condition $ d_2\alpha_{11}+d_1\alpha_{22}>2\sqrt{d_1d_2}(\alpha_{11}\alpha_{22}-\alpha_{12}\alpha_{21})$ holds, and the positivity of $d_2\alpha_{11}+d_1\alpha_{22}$ ensures the feasibility of $k_1$ and $k_2$. Depending upon the choice of parameter values and the size of the domain, we can find more than one combination of $(n_1,n_2)$ for which the following result holds
\begin{eqnarray}
k_1<\frac{\pi}{L}\sqrt{n_1^2+n_2^2}<k_2.
\end{eqnarray}
For a given set of parameter values, we can find a different admissible set of values for $(n_1,n_2)$ such that the above inequality is satisfied. But, the resulting stationary pattern is the combination of $(n_1,n_2)$  which corresponds to the maximum value of the real part of $\lambda(k)$ which is a root of the characteristic equation (\ref{cheqn1}).

\section{Stationary and non-stationary patterns}

In this work, we are interested in understanding the pattern formation over a complex domain whose shape is different from a square or a rectangle. For this purpose, we first present various types of stationary and non-stationary patterns produced by model (\ref{model1}) over a square domain briefly. A detailed discussion on these patterns can be found in our earlier works \cite{MBWWMBE,KMMBEC}. We summarize the resulting patterns produced by the model (\ref{model1}) over a square domain in the following subsection, along with the illustration of the analytical results mentioned in the previous section.

\begin{figure}[ht!]
\begin{center}
        \begin{subfigure}[p]{0.31\textwidth}
                \centering
               \includegraphics[width=\textwidth]{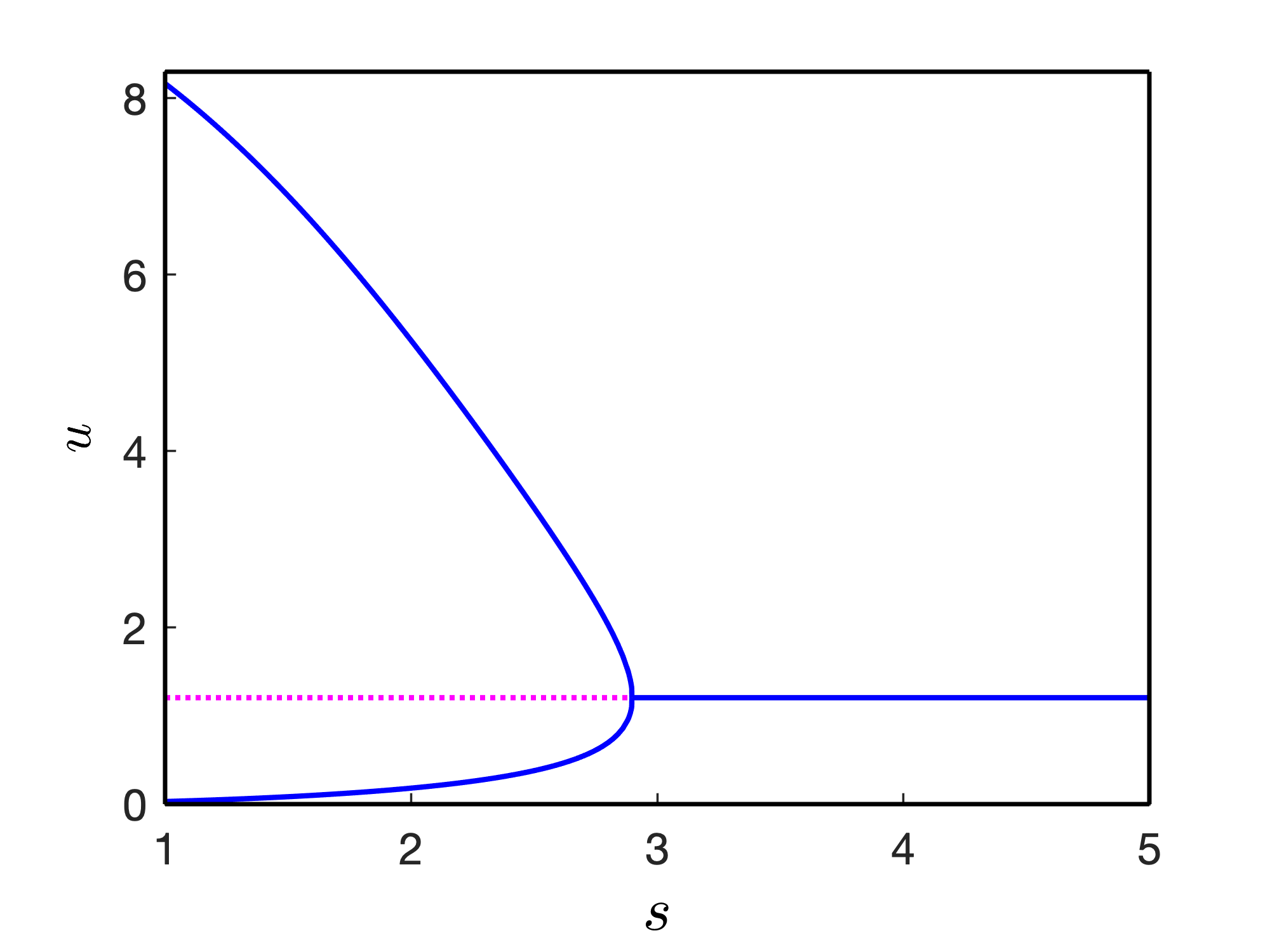}
                \caption{$a=1.5$ }\label{fig:LCa}
        \end{subfigure}%
        ~~
        \begin{subfigure}[p]{0.31\textwidth}
                \centering
               \includegraphics[width=\textwidth]{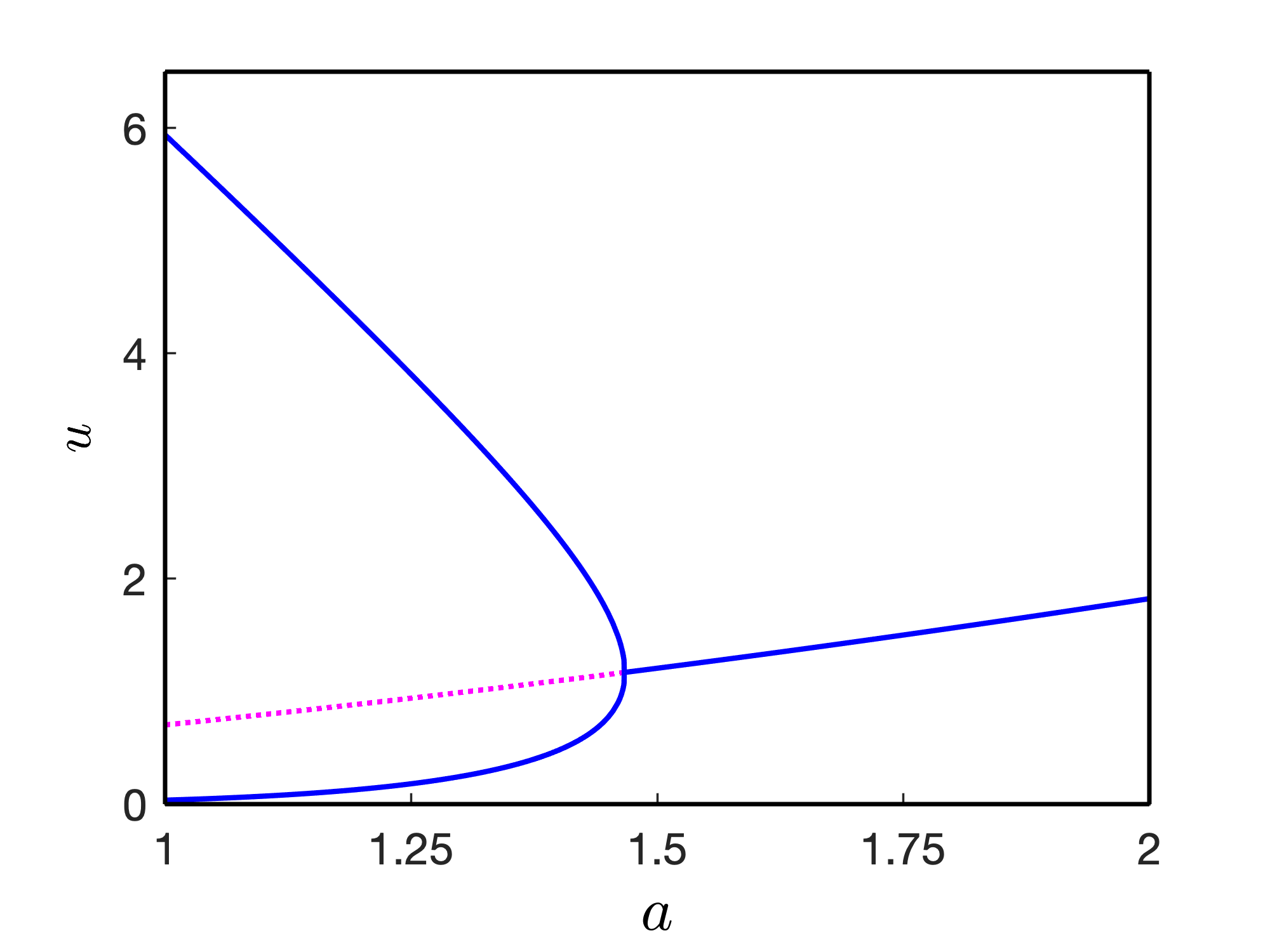}
                \caption{ $s=3$ }\label{fig:LCb}
        \end{subfigure}
        \caption{Bifurcation diagrams of the coexistence steady-state through super-critical Hopf-bifurcation with the variation of (a) the parameter $s$ and (b) the parameter $a$.}\label{fig:hopf}
\end{center}
\end{figure}

We fix the parameter values $r=1$, $f=0.1$, $m=0.1$, $b=0.9$, $c=1$, $q=0.35$, $p=0.0425$ throughout this work and vary $s$, $a$, $d_1$ and $d_2$. We vary the parameter $a$ within the range $[1.5,2.5]$ such that the model admits only one coexistence equilibrium point for the temporal model (\ref{model2}) and unique coexistence homogeneous steady-state for the spatio-temporal model (\ref{model1}). The parameter $s$ has no role to determine the level of steady-state, however, regulates the stability of coexistence state. Change in the stability of coexistence equilibrium point through super-critical Hopf-bifurcation is shown in Fig.~\ref{fig:hopf} with the variation of the parameters $a$ and $s$ respectively. We prepare two bifurcation diagrams by keeping one of the parameters $a$ and $s$ fixed while varying the other one. For $a=1.5$, the Hopf-bifurcation threshold is $s=2.89897$ and for $s=3$ fixed we find the Hopf-bifurcation threshold as $a=1.467136$. For the first case, i.e., for $a=1.5$, $s=2.89897$, we calculate the first Lyapunov number calculated as $\sigma=-0.011388$. Here, the negative first Lyapunov number ensures that the Hopf-bifurcation is super-critical.

\subsection{Pattern formation over square domain}

\begin{figure}[ht!]
\begin{center}
        \begin{subfigure}[p]{0.31\textwidth}
                \centering
                \includegraphics[width=\textwidth]{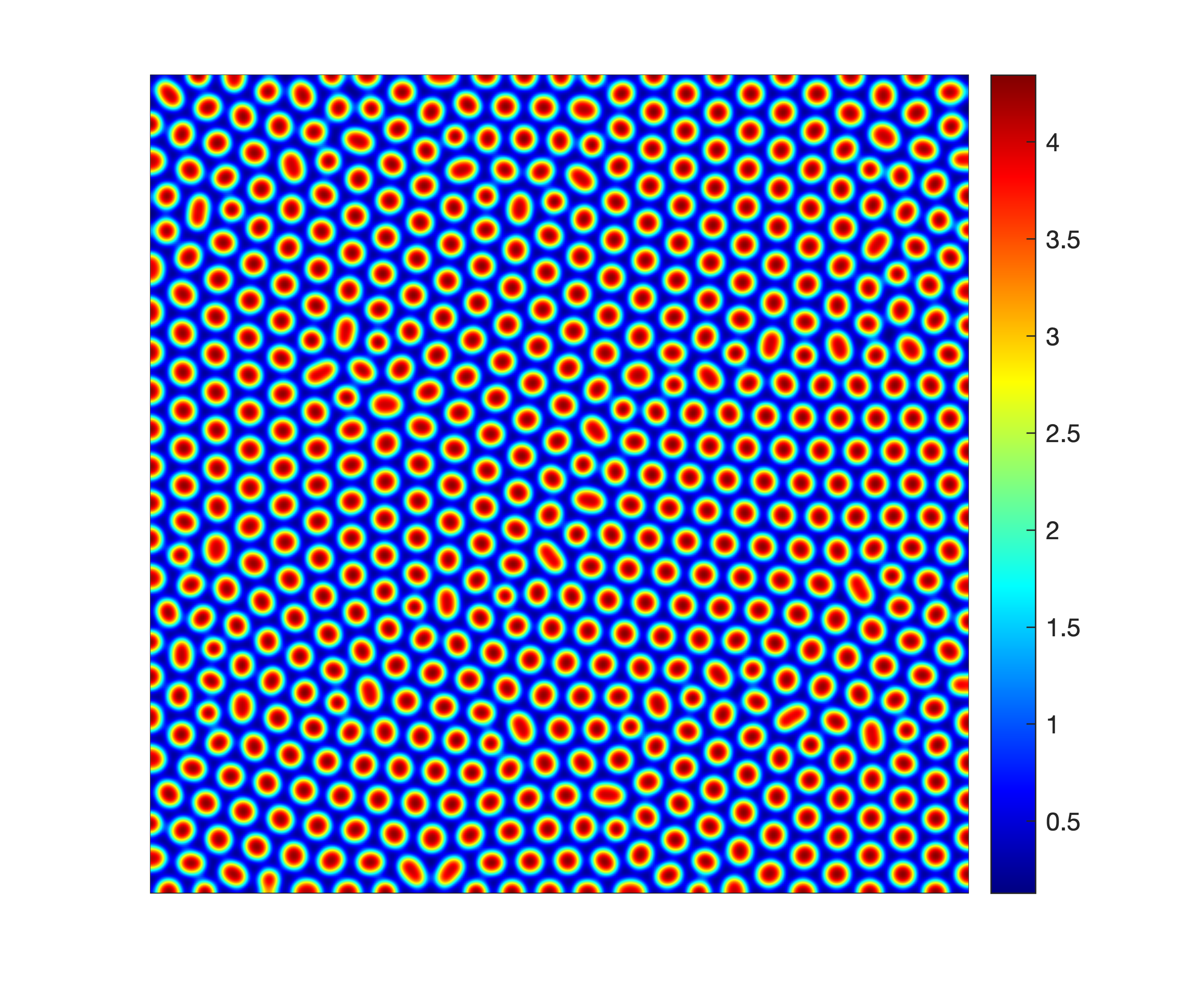}
                \caption{ }\label{fig:pat1a}
        \end{subfigure}%
        \begin{subfigure}[p]{0.31\textwidth}
                \centering
                \includegraphics[width=\textwidth]{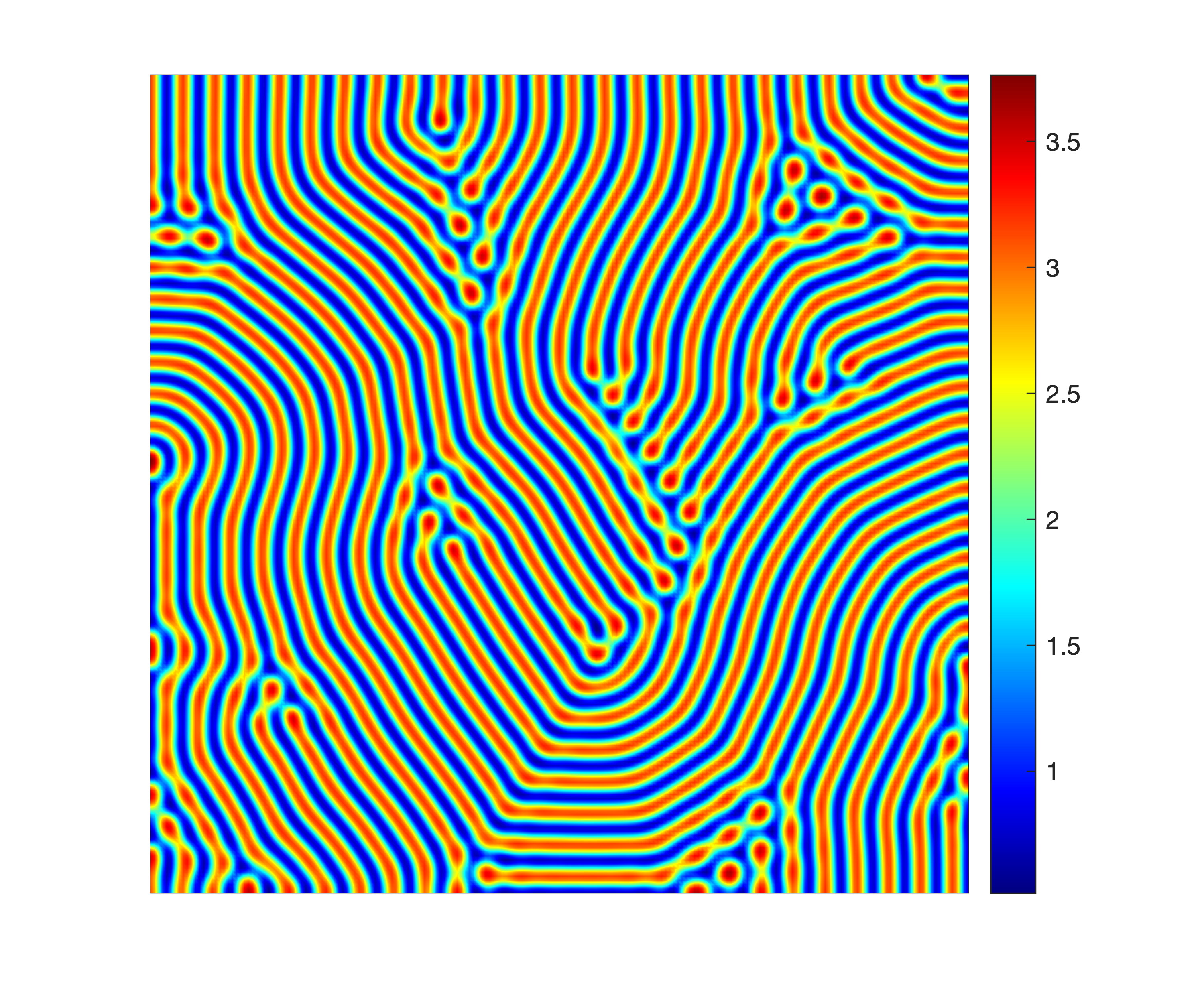}
                \caption{ }\label{fig:pat1b}
        \end{subfigure}%
        \begin{subfigure}[p]{0.31\textwidth}
                \centering
                \includegraphics[width=\textwidth]{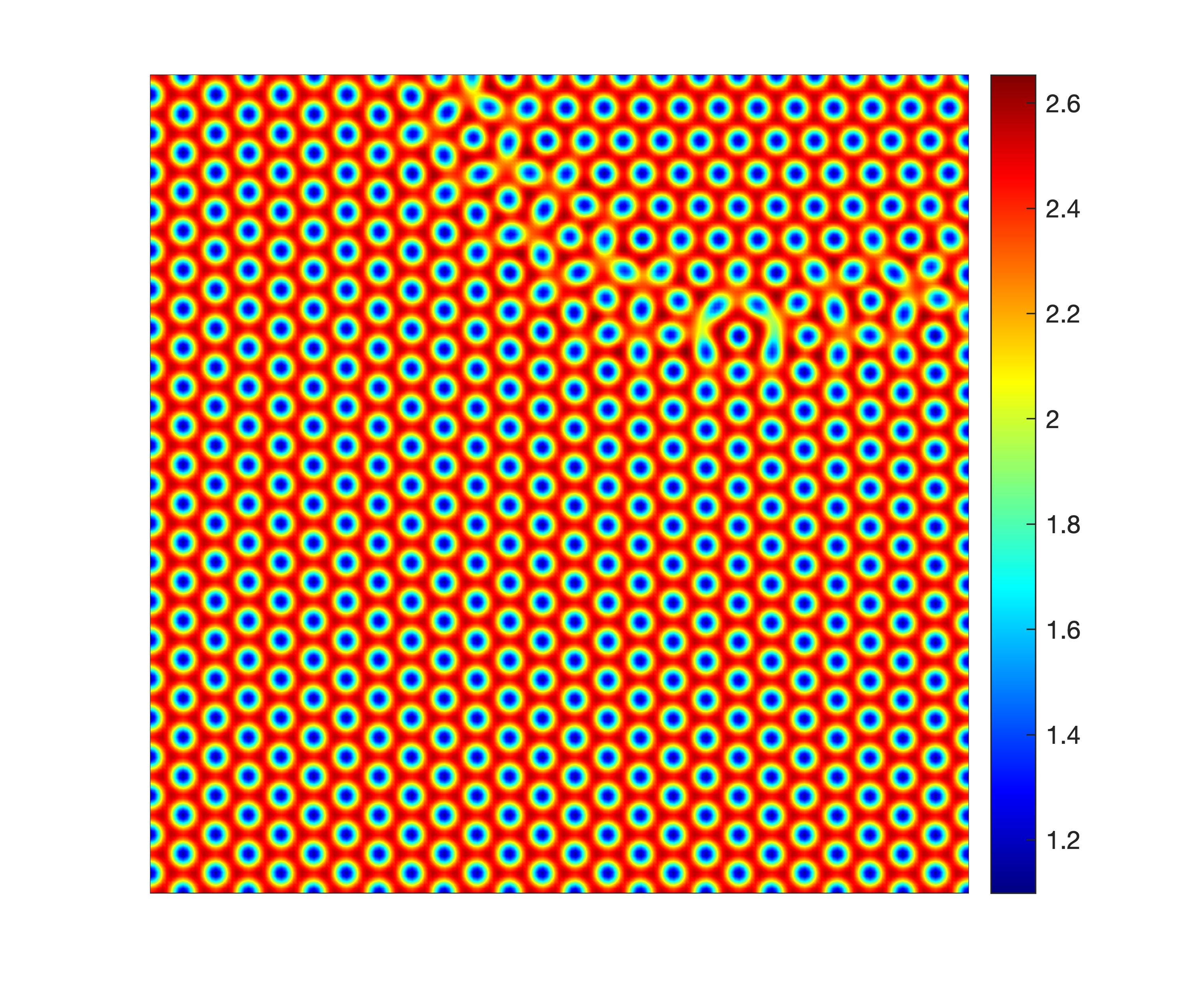}
                \caption{ }\label{fig:pat1c}
        \end{subfigure}%
\caption{Stationary Turing patterns produced by the model (\ref{model1}) for fixed parameter values $s=3$, $d_2=10$ and three different values $a$: (a) $a=1.65$, hot spot pattern; (b) $a=2$, labyrinthine pattern; (c) $a=2.2$, cold spot pattern.}\label{fig:sqdom1}
\end{center}
\end{figure}

Here we consider the pattern formation by the model (\ref{model1}) over a square domain $[0,L]\times[0,L]$, subject to zero-flux boundary condition and non-negative initial condition
\begin{eqnarray}\label{IC1_square}
u(0,x,y)=u_*+\epsilon\xi(x,y),\,\,v(0,x,y)=v_*+\epsilon\eta(x,y),
\end{eqnarray}
where $\xi(x,y)$ and $\eta(x,y)$ are two Gaussian white noises $\delta$-correlated in space \cite{HMSPEVBook, MBSPTE} and $\epsilon$ measures the amplitude of small heterogeneous perturbation from the homogeneous steady-state. As mentioned above, henceforth, we mention only the values of $s$, $a$, $d_1$ and $d_2$ to obtain various stationary and dynamic patterns. We present three different stationary patterns in Fig.~\ref{fig:sqdom1}. For stationary patterns we fix $s=3$, $d_1=0.15$, $d_2=10$, we find hot spot pattern for $a=1.65$, labyrinthine pattern for $a=2$ and cold spot pattern for $a=2.2$. We choose the length of the square domain as $L=200$ for the numerical solution of stationary patterns. It is important to mention here that the obtaining stationary patterns presented in Fig.~\ref{fig:sqdom1} are at a different time as the transient pattern persists for different duration of time. Three different stationary patterns are obtained at $t=500$ (hot spot), $t=1800$ (labyrinthine) and $t=2500$ (cold spot) respectively. For these choices of parameter values, homogeneous coexistence steady-state is stable, and the magnitude of diffusion coefficients satisfy the Turing instability condition. Hence, the resulting patterns are Turing patterns and are stationary. However, the duration of the transient pattern depends on the strength of temporal interactions and the magnitude of the diffusion parameters.

\begin{figure}[ht!]
\begin{center}
        \begin{subfigure}[p]{0.31\textwidth}
                \centering
                \includegraphics[width=\textwidth]{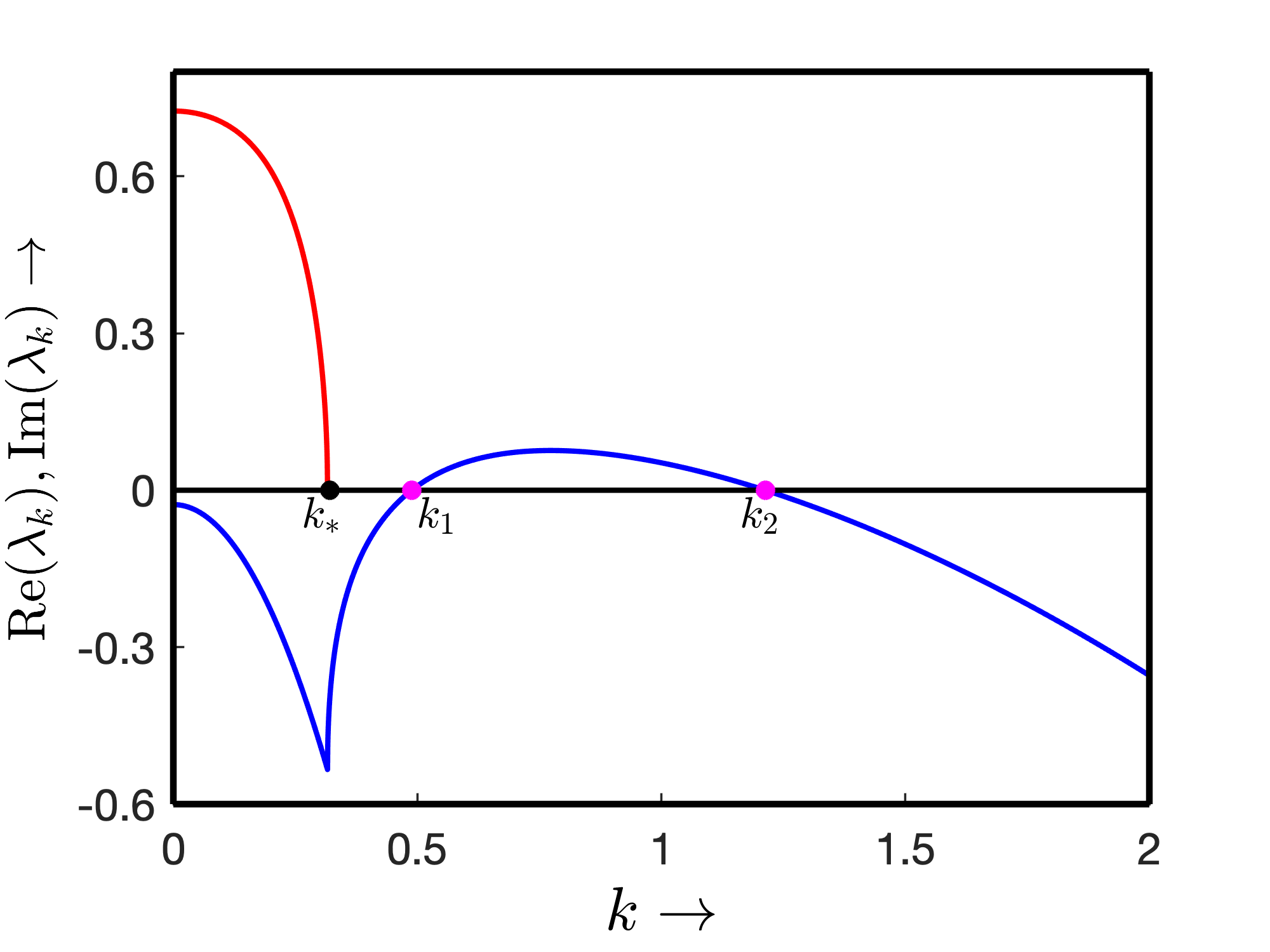}
                \caption{ }\label{fig:dispa}
        \end{subfigure}%
        \begin{subfigure}[p]{0.31\textwidth}
                \centering
                \includegraphics[width=\textwidth]{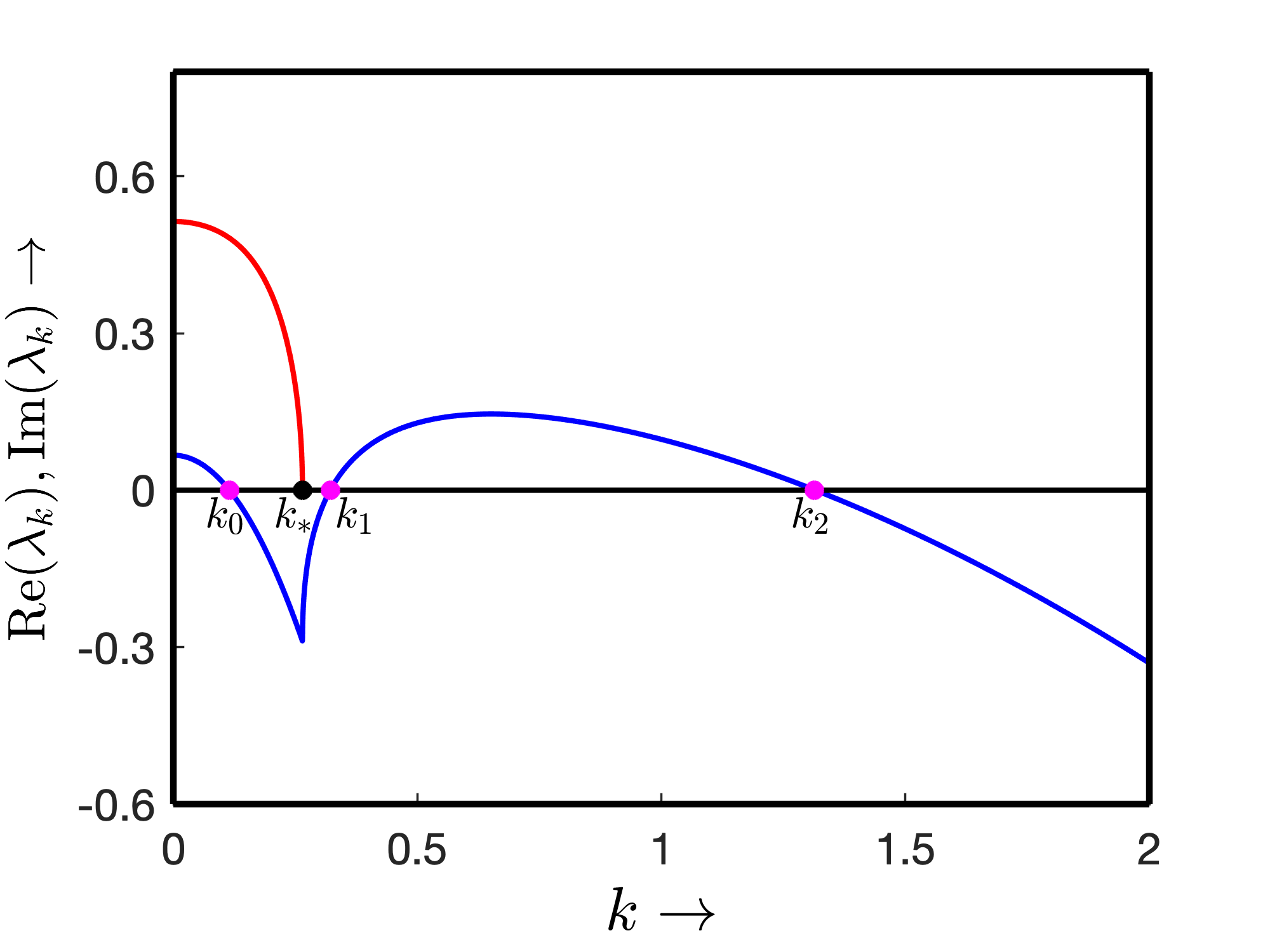}
                \caption{ }\label{fig:dispb}
        \end{subfigure}%
        \begin{subfigure}[p]{0.31\textwidth}
                \centering
                \includegraphics[width=\textwidth]{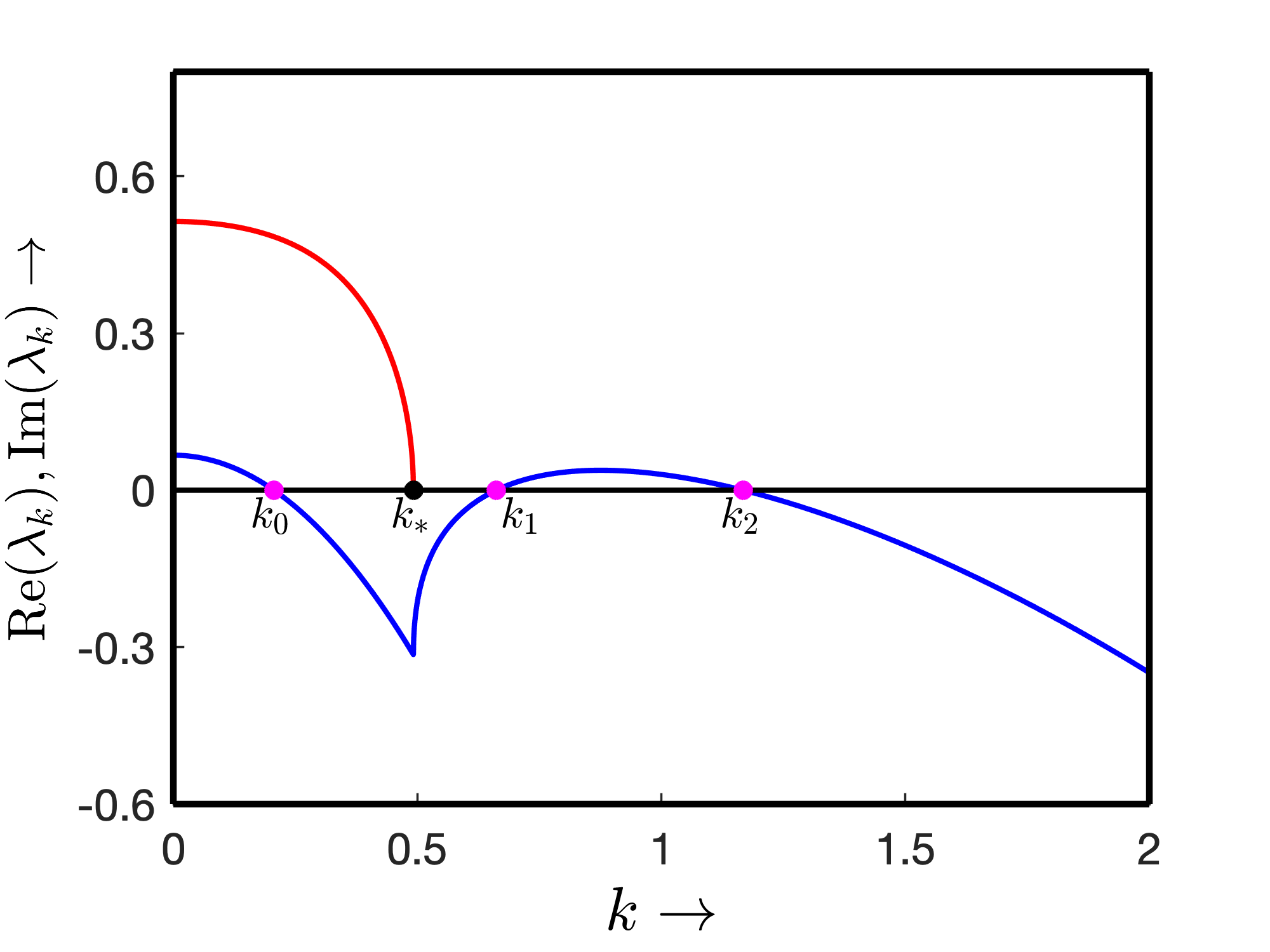}
                \caption{ }\label{fig:dispc}
        \end{subfigure}%
\caption{The plot of dispersion relations corresponding to three stationary patterns presented in Fig.~\ref{fig:sqdom1}. Imaginary parts are plotted in red and the maximum real parts plotted in blue for the eigenvalues $\lambda(k)$.}\label{fig:disp}
\end{center}
\end{figure}

The dispersion relations corresponding to three stationary Turing patterns presented in Fig.~\ref{fig:sqdom1}, are presented in Fig.~\ref{fig:disp}. For $a=1.65$, the largest real part of the eigenvalues of the characteristic equation (\ref{cheqn1}) is positive for $k\in(k_1,k_2)=(0.49,1.21)$ and is maximum at $k=0.773$. There are plenty of choices for $n_1$ and $n_2$ such that the inequality $0.49<\frac{\pi}{L}\sqrt{n_1^2+n_2^2}<1.21$ holds, for $L=200$. In particular, the value $k=0.773$, corresponding to the maximum real part of the eigenvalue, can be achieved from the expression $\frac{\pi}{L}\sqrt{n_1^2+n_2^2}$ for a couple of choices of $n_1$, $n_2$, namely $(n_1,n_2)=(35,35)$, $(n_1,n_2)=(34,36)$, $(n_1,n_2)=(36,34)$ and so on. The choice of $(n_1,n_2)$ for which $\frac{\pi}{200}\sqrt{n_1^2+n_2^2}$ is close to $k=0.773$ is $(n_1,n_2)=(35,35)$. The resulting pattern presented in Fig.~\ref{fig:sqdom1}(\subref{fig:pat1a}) shows equal number of peaks within the domain $[200\times200]$ as the number of peaks we find for the function $\cos\frac{n_1\pi x}{200}\cos\frac{n_2\pi x}{200}$ when plotted over a same domain size. Similar argument holds for the patterns presented in other panels of Fig.~\ref{fig:sqdom1}.

\begin{figure}[ht!]
        \centering
        \includegraphics[width=110mm]{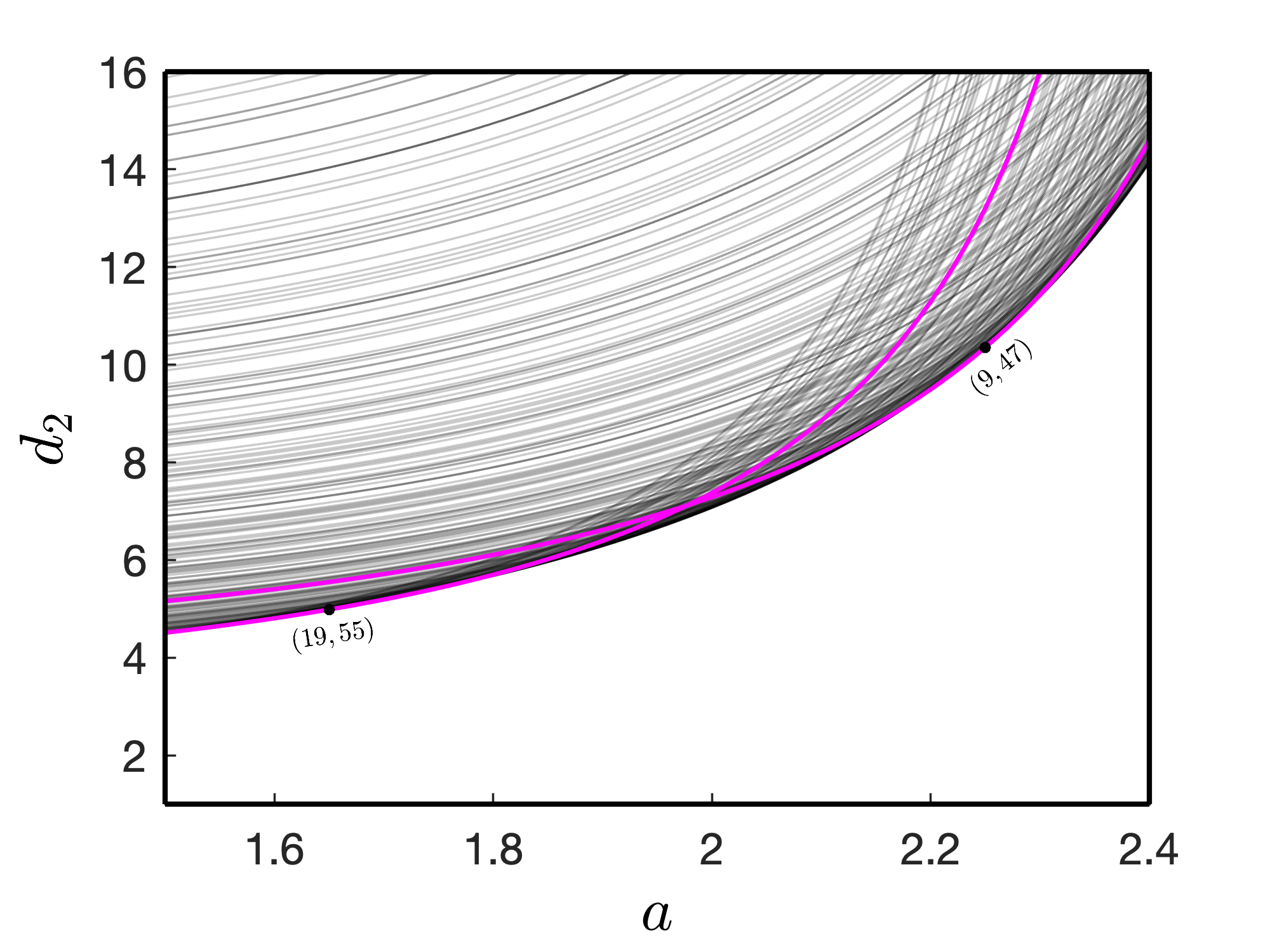}
\caption{Critical Turing instability curves for various unstable spatial modes $(n_{1},n_{2})$ according to the stability condition $h(k^{2})>0$ defined in (\ref{eq:h_k}). The fixed parameter values are $r=1$, $f=0.1$, $m=0.1$, $b=0.9$, $c=1$, $s=3$, $q=0.35$, $p=0.0425$ and $d_1=0.15$.}\label{fig:kmn}
\end{figure}

We present some critical curves for the Turing instability of the spatially homogeneous steady state in Fig.~\ref{fig:kmn}. We plot these curves in thin black. For each $(n_{1},n_{2})$, we find these curves by solving the equation (\ref{eq:h_k}) for $d_{2}$ with varying $a$. Specifically, two critical curves are marked (thick magenta curves) with the mode numbers, while the rest are not (thin black curves). For $a=1.65$, the Turing bifurcation threshold lies close to the critical curve mode $(19,55)$. It is the only exciting mode that determines the spatial pattern of the emerging wave dynamics. For $a=1.65$, $d_{2}=4.986$, the largest real part of the eigenvalues of the characteristic equation (\ref{cheqn1}) is positive for $k\in(k_1,k_2)=(0.911,0.921)$ and is maximum at $k=0.916$. Therefore, for $n_1=19$, $n_2=55$ and $L=200$, the inequality $0.911<\frac{\pi}{L}\sqrt{n_1^2+n_2^2}<0.921$ holds as $k=0.916$ corresponds to the maximum real part of the eigenvalue, and it is close to the value of the expression $\frac{\pi}{L}\sqrt{n_1^2+n_2^2}=0.914$ while $(n_1,n_2)=(19,55)$. Hence, $(19,55)$ is a reasonable critical mode and the corresponding mode curve is close to the Turing bifurcation threshold \cite{TZA}. Similarly, the mode curve for $a=2.25$ and close to Turing bifurcation curve is determined by the mode numbers $(n_1,n_2)=(9,47)$.

\begin{figure}[ht!]
        \centering
        \includegraphics[width=150mm]{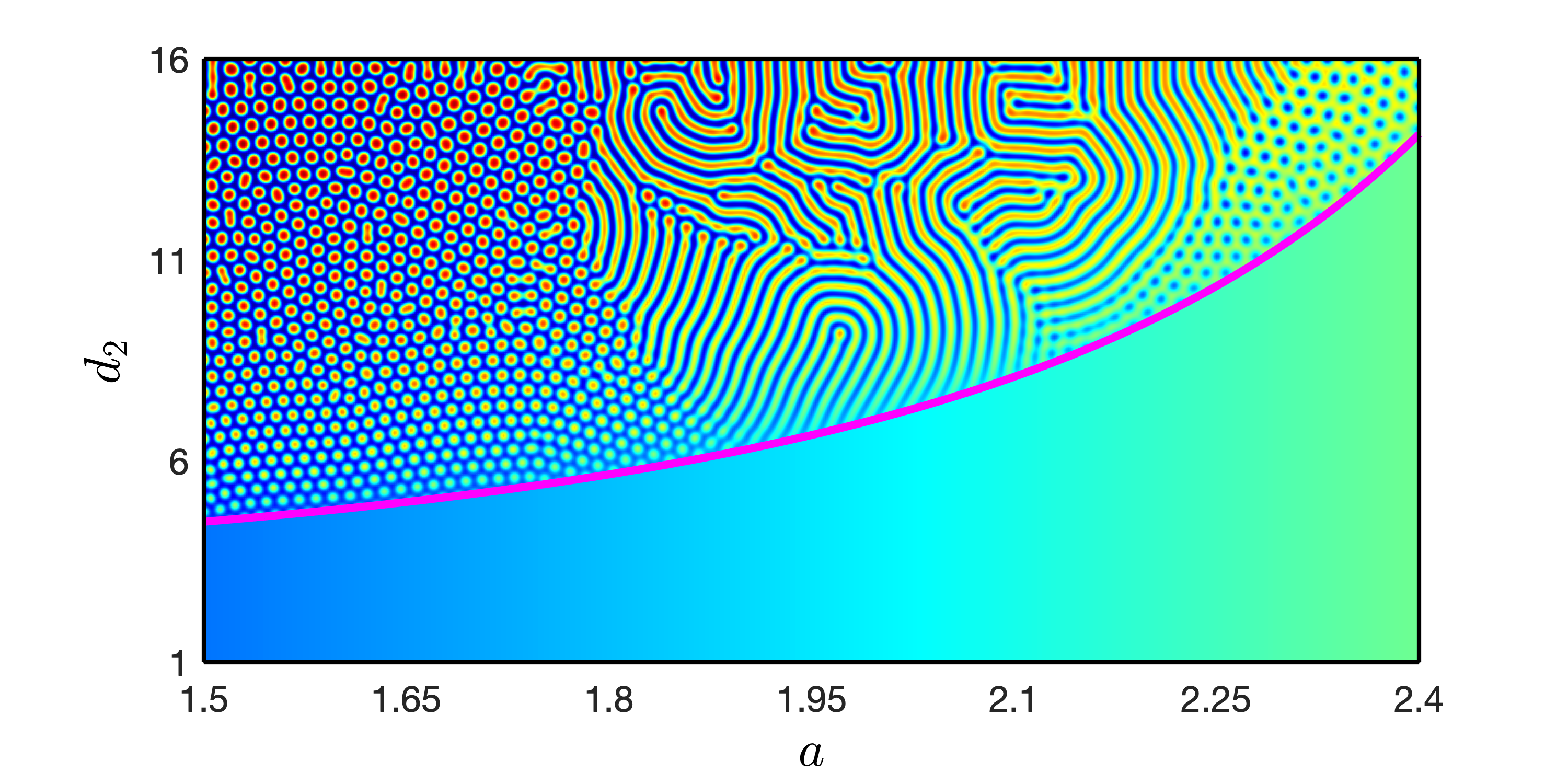}
\caption{Consolidated pattern diagram in $a-d_2$-parameter space of the spatio-temporal model (\ref{model1}) for the parameter values $s=3$ and $d_1=0.15$. }\label{fig:SchematicPattern}
\end{figure}

The stationary patterns presented in Fig.~\ref{fig:sqdom1} correspond to the parameter values from the pure Turing domain as well as from Turing-Hopf domain. From the dispersion relations, we can identify that $a=1.65$ corresponds to the parameter value in the pure Turing domain as the real part of $\lambda(k)$ is negative at $k=0$. On the other hand, for $a=2$ and $a=2.5$, the real parts of $\lambda(k)$ is positive when $k=0$. For parameter values within the Turing-Hopf domain, the real part of the eigenvalue for $k=0$ is less than the maximum real part of the eigenvalue for some $k \neq 0$, and hence we find stationary patterns. A consolidated pattern diagram is presented in Fig.~\ref{fig:SchematicPattern} for a range of values of the parameters $a$ and $d_2$. Stationary heterogeneous patterns exist for parameter values above the Turing bifurcation curve (magenta color curve). Below the Turing bifurcation curve, homogeneous steady-state is stable, and the level of stable homogeneous steady-state varies with the change of the value of $a$ from 1.5 to 2.4.

\begin{figure}[ht!]
\begin{center}
        \begin{subfigure}[p]{0.35\textwidth}
                \centering
                \includegraphics[width=\textwidth]{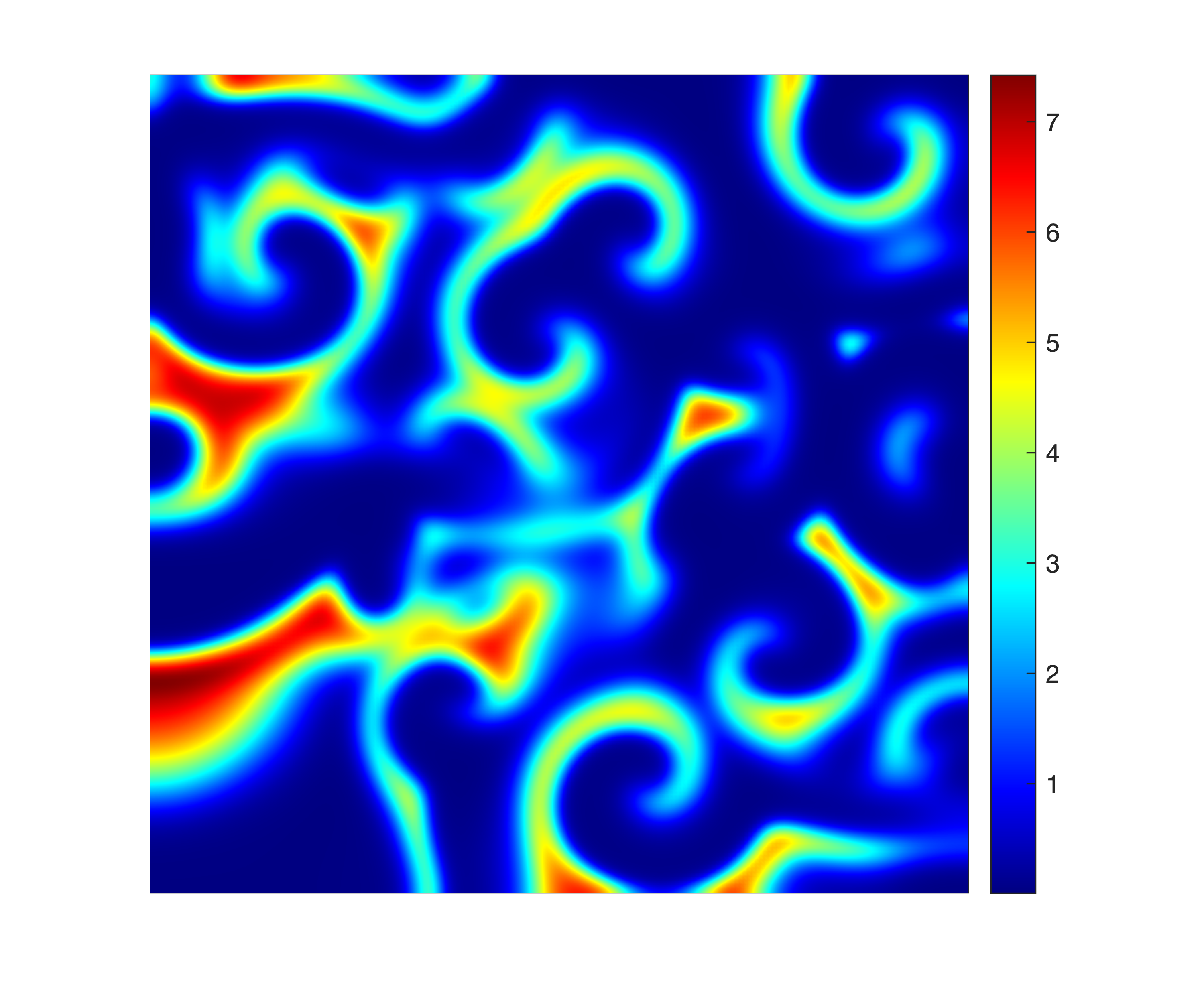}
                \caption{ }
        \end{subfigure}%
        \begin{subfigure}[p]{0.35\textwidth}
                \centering
                \includegraphics[width=\textwidth]{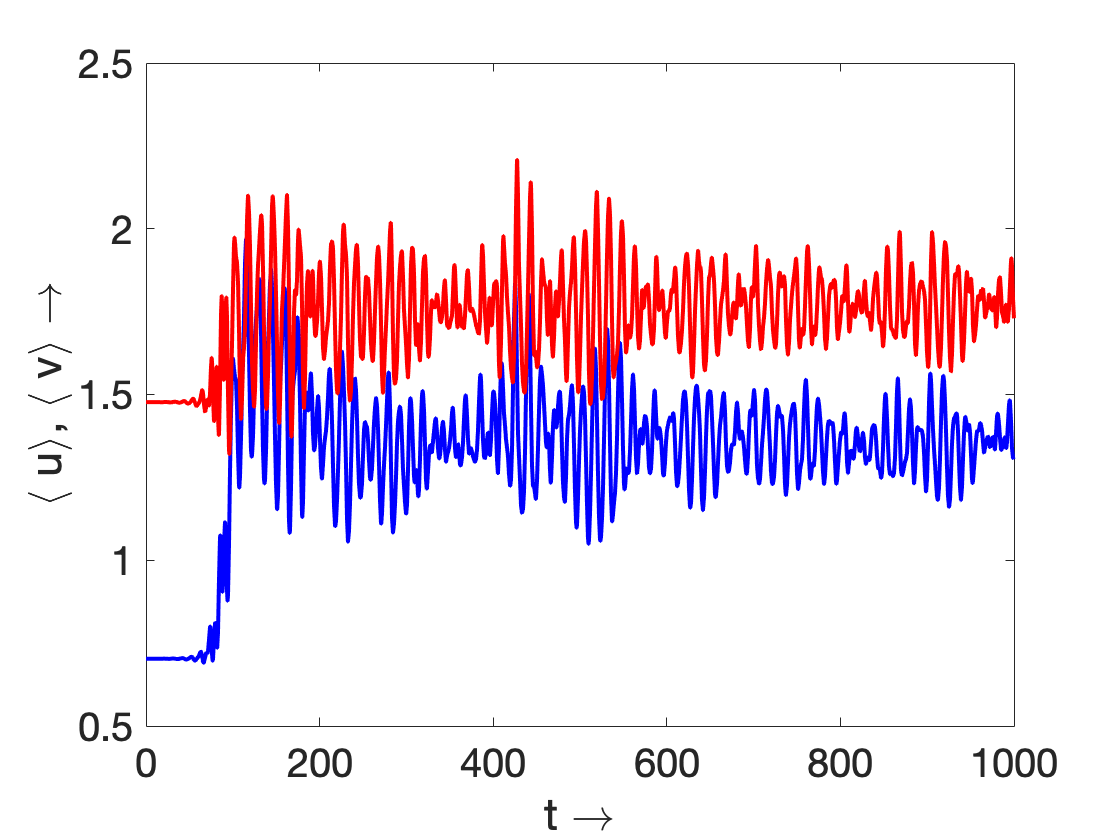}
                \caption{ }
        \end{subfigure}
\caption{(a) Spatio-temporal chaotic pattern produced by the model (\ref{model1}) for the parameter values  $a=1.0$, $s=1$, $d_1=1$ and $d_2=1$. (b) Time evolution of spatial average of prey (blue) and predator (red) densities for $t\in[0,1000]$.}\label{fig:sqdom2}
\end{center}
\end{figure}

Now we consider a non-Turing pattern for $d_1=d_2=1$ which is in fact a non-stationary pattern. For $a=1.0$, $s=2.0$, and other parameter values as mentioned earlier, we find spatio-temporal chaos as shown in Fig.~\ref{fig:sqdom2}(a). We obtain this spatio-temporal chaotic pattern through the numerical simulation over a square domain of size $L=200$ by taking the initial condition as mentioned in (\ref{IC1_square}) and the zero-flux boundary condition. The chaotic pattern is shown for $t=1000$, the nature of the chaotic pattern is interacting spiral chaos. We have verified the chaotic nature of the resulting pattern by plotting the chaotic oscillation of time evolution of spatial averages of prey and predator population densities $\langle u\rangle$ and $\langle v\rangle$ in Fig.~\ref{fig:sqdom2}(b). This spatio-temporal chaotic pattern also exhibits sensitivity to the initial condition, but we omit that result for the sake of brevity. A detailed discussion about the verification of sensitivity of spatial pattern to initial condition for the spatio-temporal model is available in \cite{MBSAEC}.

\subsection{Pattern formation in U-shaped domain}

\begin{figure}[ht!]
        \centering
        \includegraphics[width=90mm]{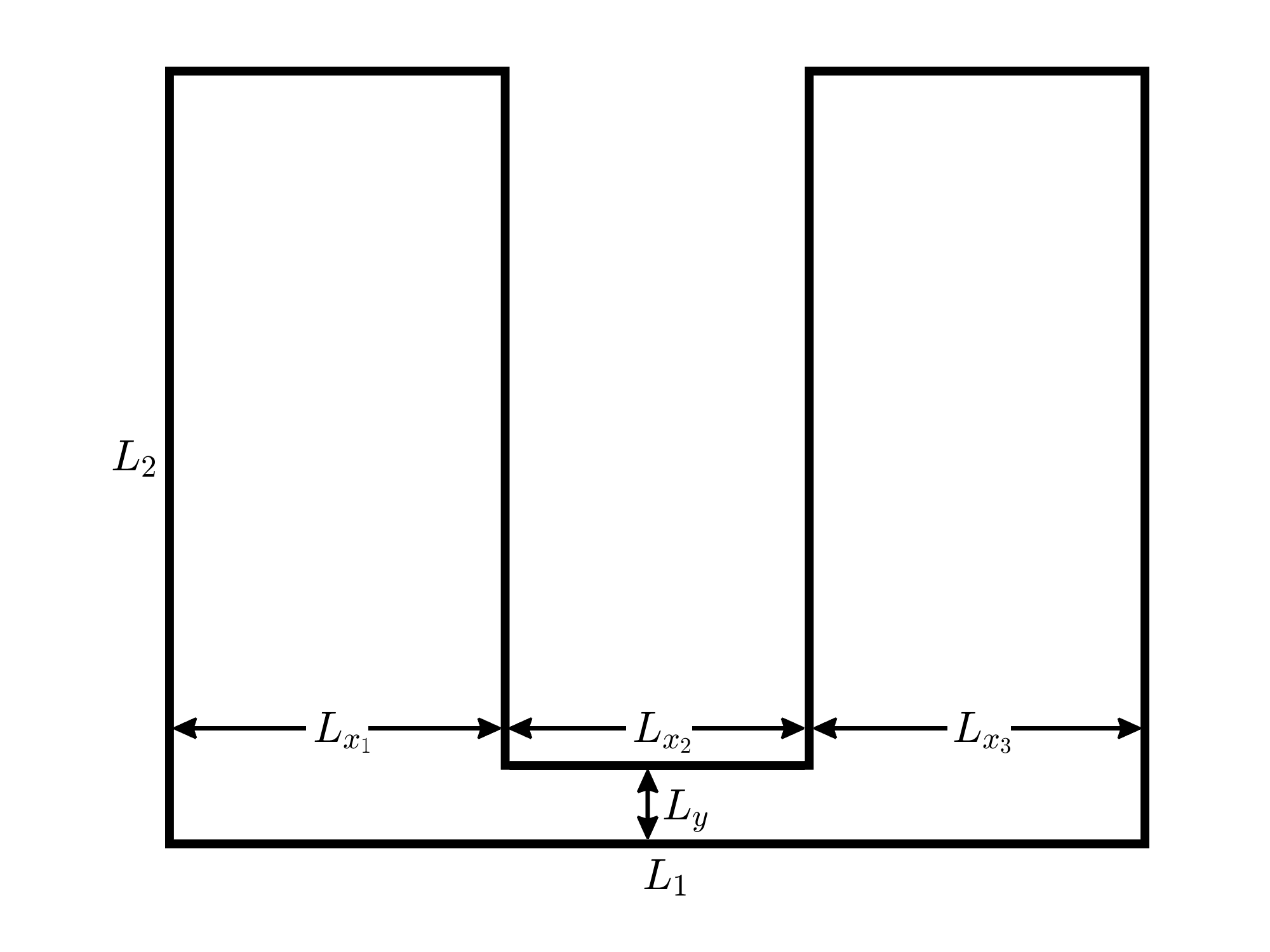}
\caption{Illustration for the U-shaped spatial domain.}\label{fig:domain}
\end{figure}

In this section, we consider the effect of the shape of the domain on spatio-temporal pattern formation. To understand the effect of the complex shape of the domain on the pattern formation, we consider the U-shaped domain as shown in Fig.~\ref{fig:domain}. Here, the connection between two large habitats is like a channel where the individuals can move from one domain to another. The U-shaped domain is characterized by the length ($L_2$) and width ($L_{x_1}$, $L_{x_3}$) of two habitats and the length ($L_{x_2}$) and width ($L_y$) of the connecting channel. First, we consider the pattern formation over the U-shaped domain when the size of the two patches are the same, and the connecting channel is reasonably wide. Next, we consider how the length and width of the connecting channel and the difference in the size of the two patches affect the various aspects of spatio-temporal pattern formation.

For the numerical simulations, we consider the following initial condition:
\begin{eqnarray}\label{IC1}
u(0,x,y)=\left\{\begin{array}{ll}
u_*+\epsilon\xi(x,y) &  (x,y)\in\mathcal{D}_1,\\
u_* & \textrm{elsewhere} \\
\end{array}\right.\mbox{and}~
v(0,x,y)=\left\{\begin{array}{ll}
v_*+\epsilon\eta(x,y) & (x,y)\in\mathcal{D}_1,\\
v_* & \textrm{elsewhere} \\
\end{array}\right.,
\end{eqnarray}
where $\mathcal{D}_1=[0, L_{x_1}]\times[0, L_2]$ and $\xi(x,y)$, $\eta(x,y)$ are spatially uncorrelated Gaussian white noise terms. The boundary condition is assumed to be no-flux along the boundaries of the entire domain. For our convenience of further discussion, we denote
the domain consisting of the connecting corridor and the patch on the right-hand side as $\mathcal{D}_2$. It can be defined as $\mathcal{D}_2=\{[L_{x_1}, L_{x_1}+L_{x_2}]\times[0, L_y]\}\cup\{[L_{x_1}+L_{x_2}, L_1]\times[0, L_2]\}$.

\begin{figure}[htpp]
\begin{center}
        \begin{subfigure}[p]{0.31\textwidth}
                \centering
                \includegraphics[width=\textwidth]{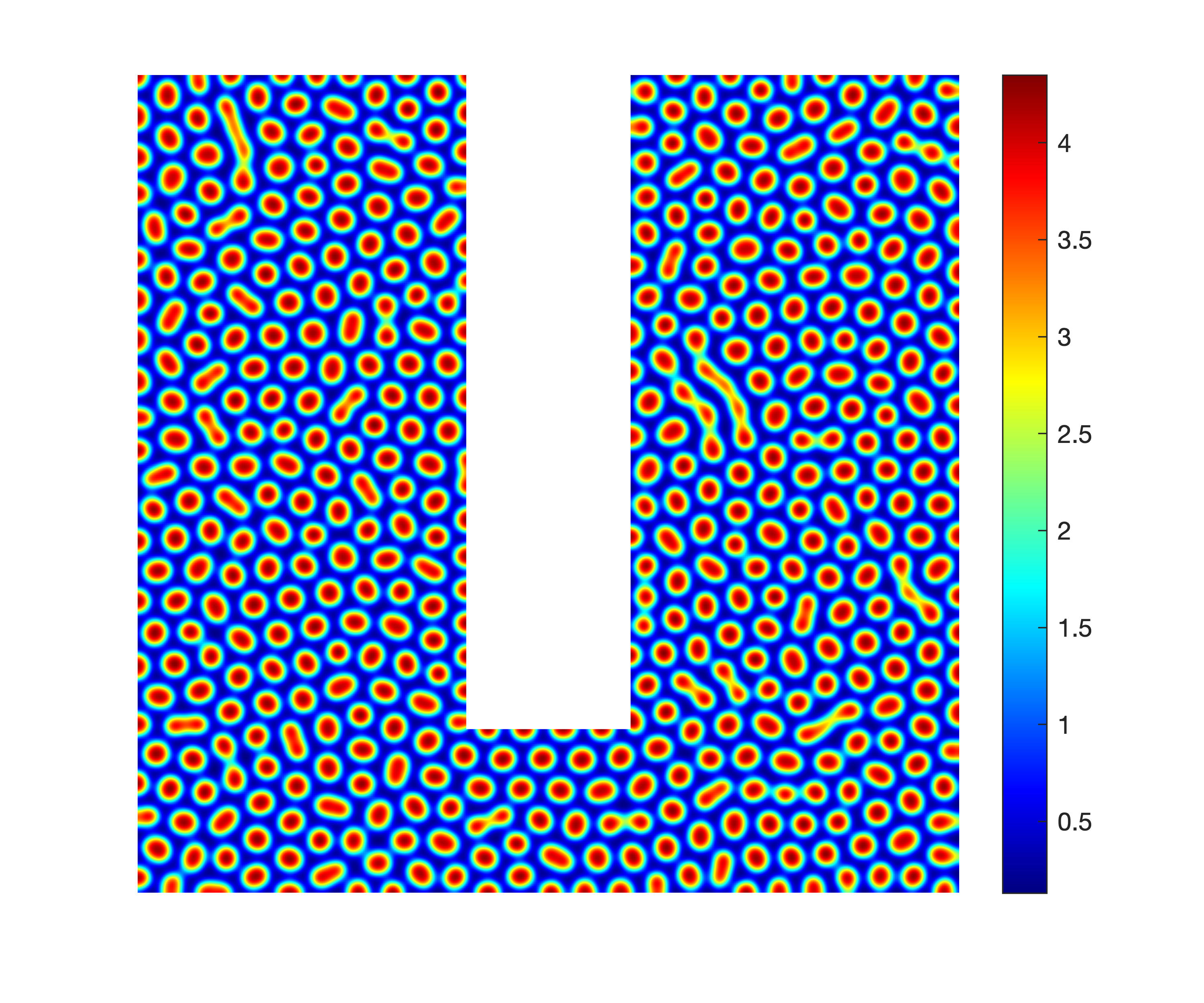}
                \caption{ }\label{fig:pat1a_sp_U}
        \end{subfigure}%
        ~~
        \begin{subfigure}[p]{0.31\textwidth}
                \centering
                \includegraphics[width=\textwidth]{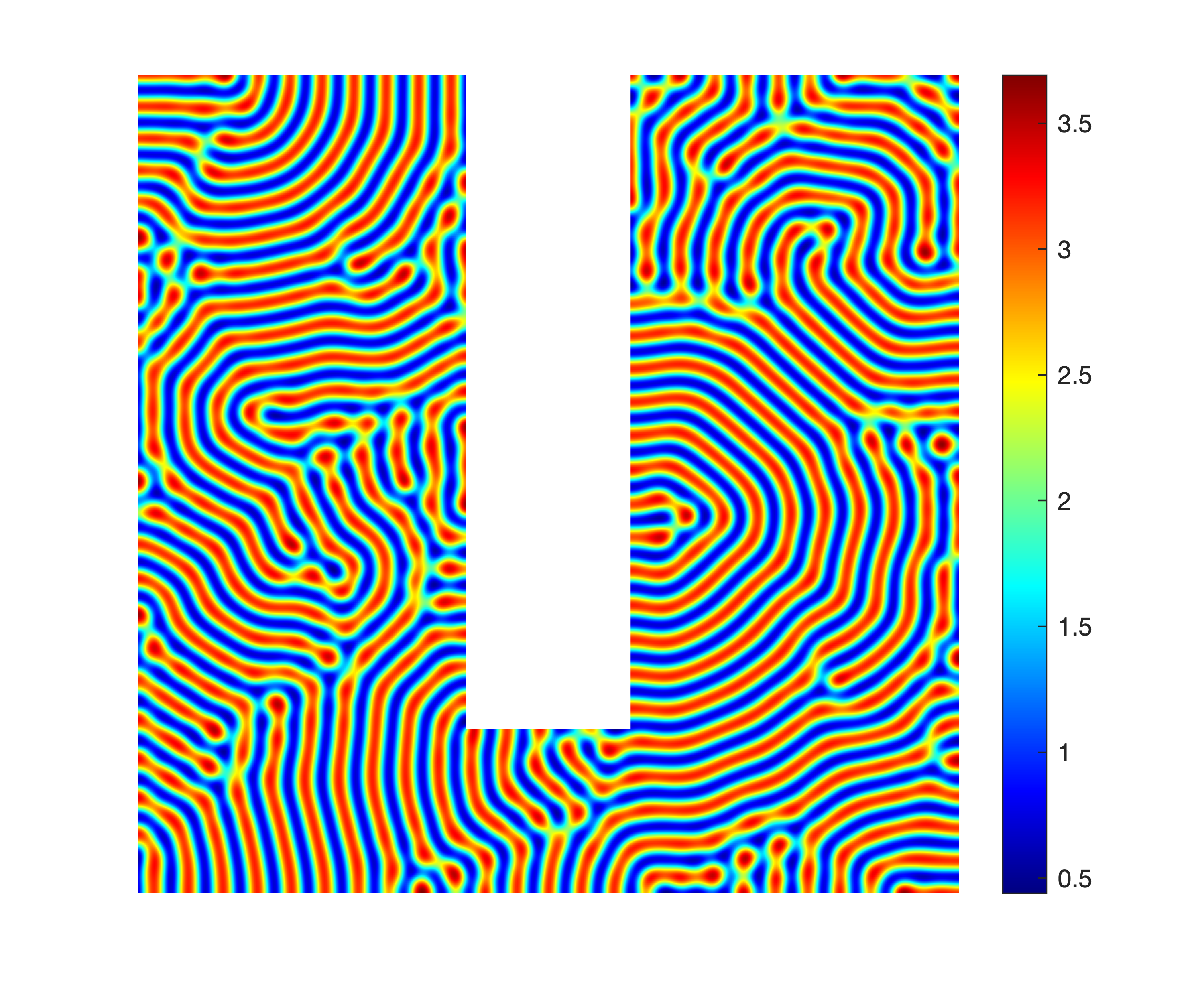}
                \caption{ }\label{fig:pat1b_sp_U}
        \end{subfigure}
        ~~
        \begin{subfigure}[p]{0.31\textwidth}
                \centering
                \includegraphics[width=\textwidth]{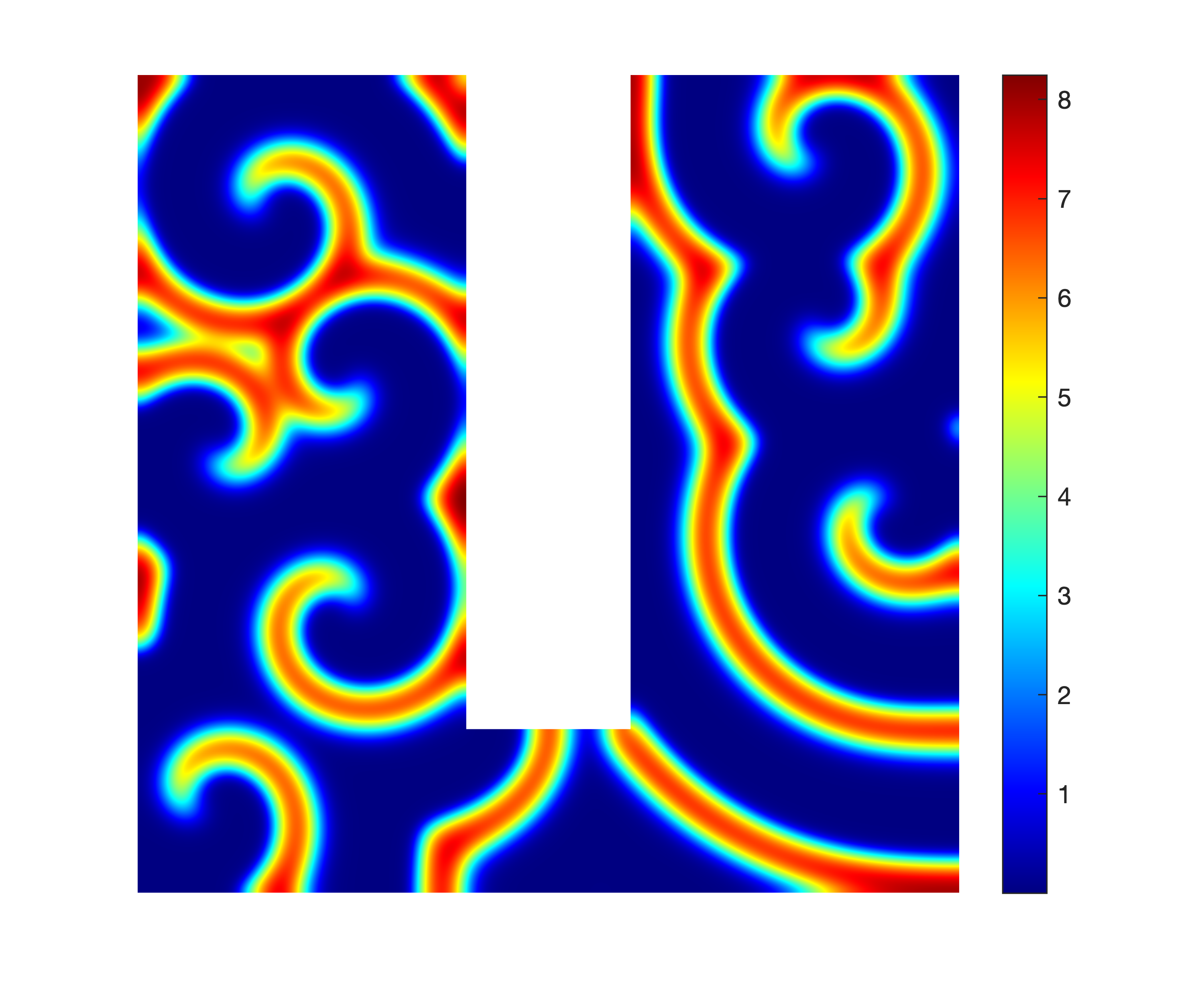}
                \caption{ }\label{fig:pat1c_sp_U}
        \end{subfigure}
        \caption{Spatio-temporal patterns over U-shaped domain for the fixed parameter values: (a) hot spot pattern ($s=3$, $a=1.5$, $d_1=0.15$, $d_2=10$); (b) labyrinthine pattern ($s=3$, $a=2$, $d_1=0.15$, $d_2=10$); (c) spatio-temporal chaotic pattern ($s=1$, $a=1$, $d_1=1$, $d_2=1$). Patterns are obtained at $t=1500$, see the text for the measure of the U-shaped domain.}\label{fig:Ushape1}
\end{center}
\end{figure}

Here, we first consider both the patches identical, and the parameter values are uniform over the entire U-shaped domain. Figure \ref{fig:Ushape1} depicts the pattern formation over the entire domain for different values of $a$, $s$, $d_1$ and $d_2$ keeping other parameters fixed as mentioned in the previous sub-section. We see that a similar type of stationary or dynamic pattern covers the entire U-shaped domain for fixed parameter values, as if the domain shape does not affect the resulting pattern. This observation is true if we perform the simulation over a reasonably long time and the width of the patches and the connecting channel is reasonably large. The characterization of the phrase `reasonably large' is explained below in detail for the specific choices of $L_{x_j}$, $j=1,2,3$ and $L_y$. The domain size in Fig.~\ref{fig:Ushape1} is $L_2=200$, $L_{x_1}=L_{x_3}=80$, $L_{x_2}=40$ and $L_y=40$. Depending on the choice of the parameters, the stationary and dynamic pattern engulf the entire U-shaped domain for the case of a reasonable large connecting channel, like the shape of the domain does not affect the pattern formation.

\begin{figure}[ht!]
\begin{center}
        \begin{subfigure}[p]{0.31\textwidth}
                \centering
                \includegraphics[width=\textwidth]{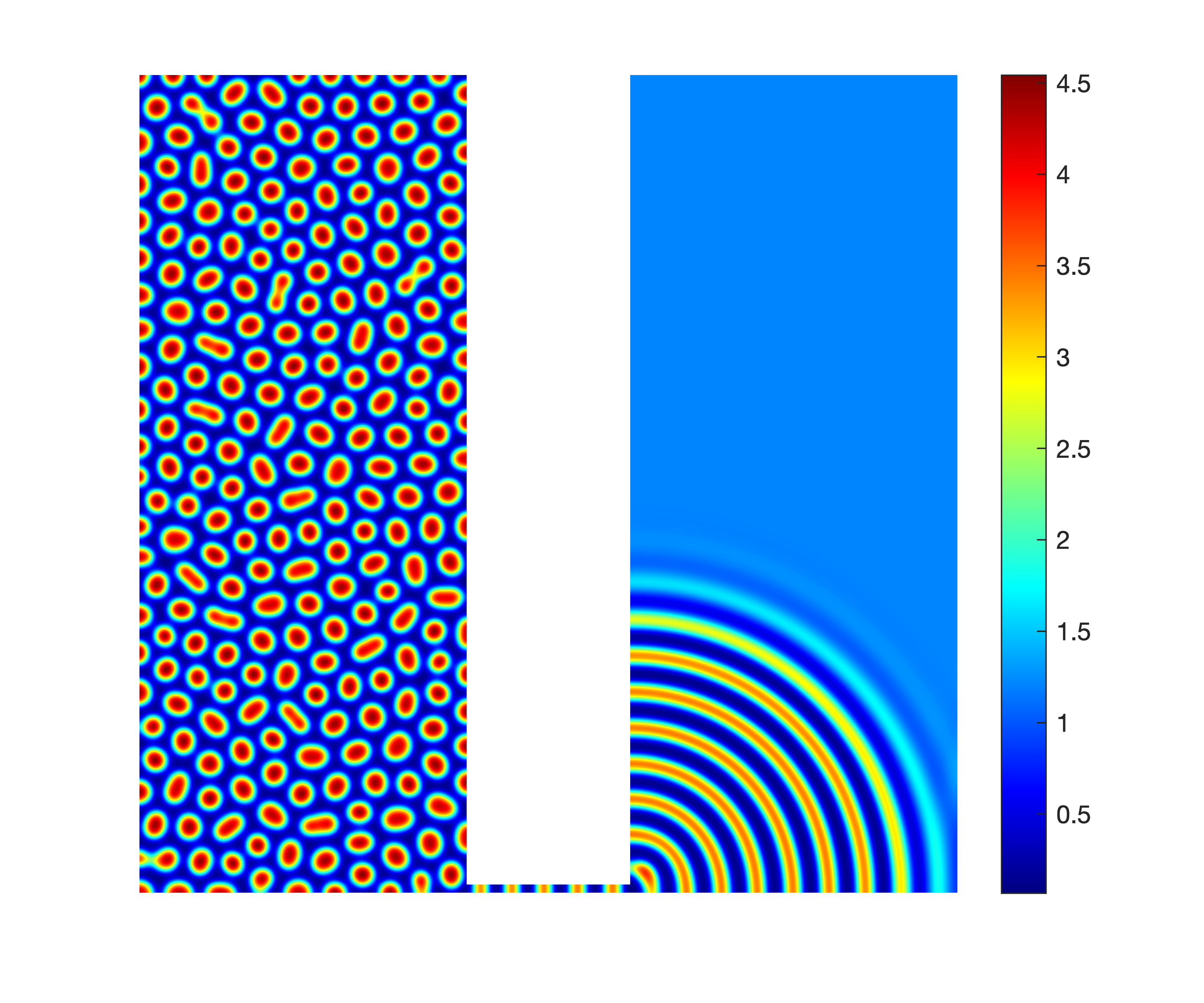}
                \caption{ }\label{fig:pat1a_sp_hotspot}
        \end{subfigure}%
        ~~
        \begin{subfigure}[p]{0.31\textwidth}
                \centering
                \includegraphics[width=\textwidth]{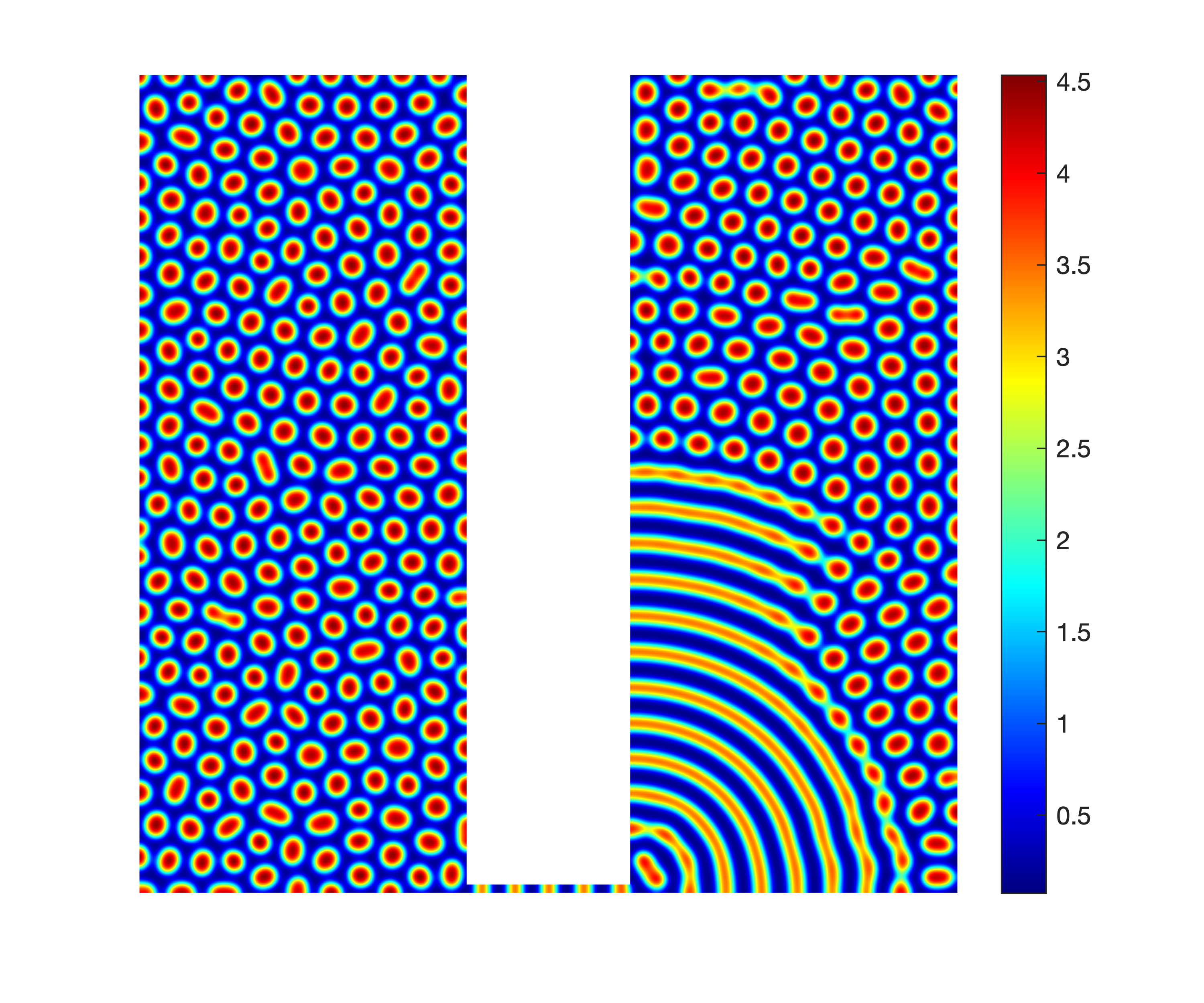}
                \caption{ }\label{fig:pat1b_sp_hotspot}
        \end{subfigure}
        \begin{subfigure}[p]{0.31\textwidth}
                \centering
                \includegraphics[width=\textwidth]{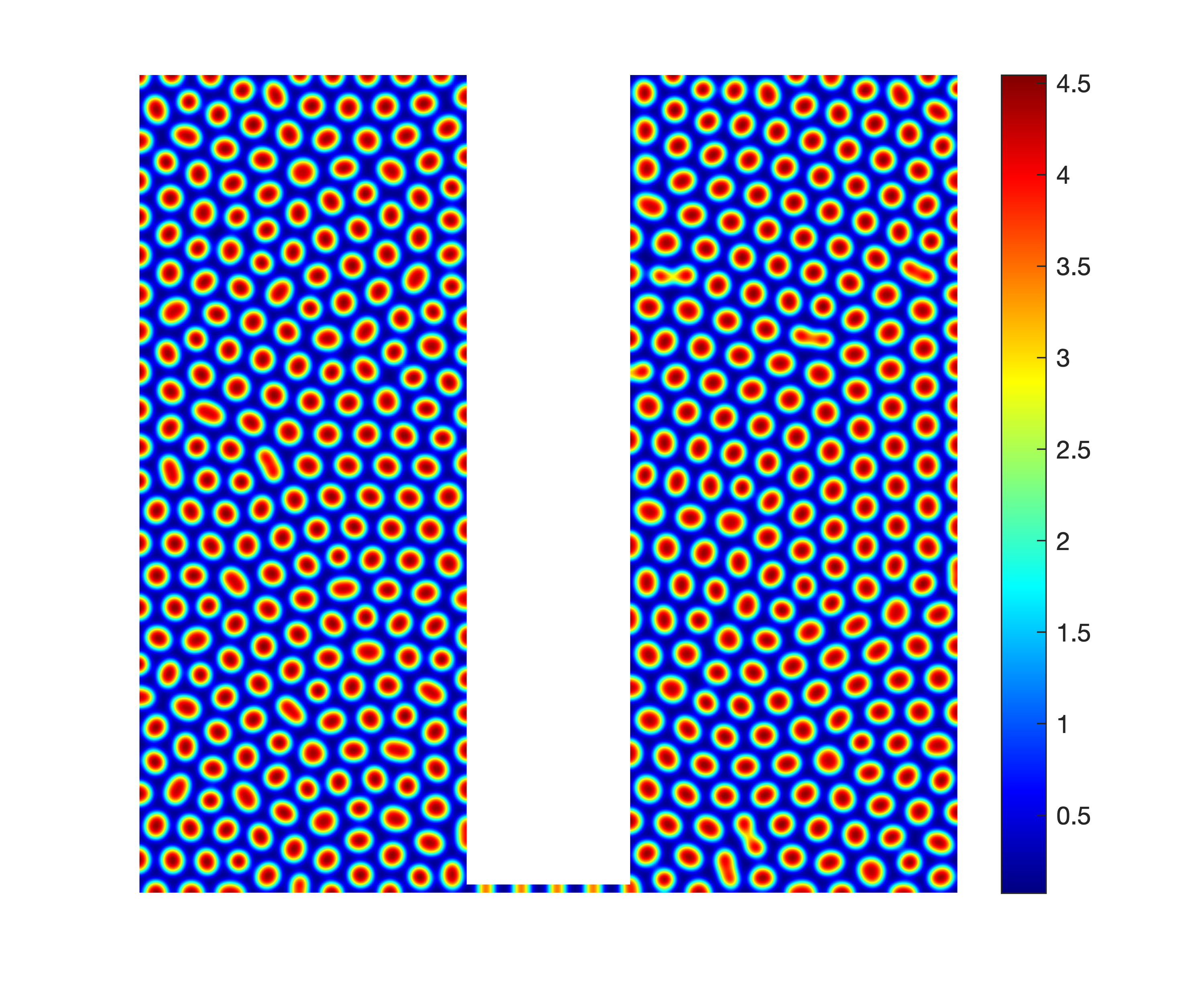}
                \caption{ }\label{fig:pat1c_sp_hotspot}
        \end{subfigure}%
        \caption{Stationary hot-spot pattern formation in U-shaped domain for the parameter values $s=3$, $a=1.5$, $d_1=0.15$ and $d_2=10$. Patterns are obtained at different time (a) $t=300$, (b) $t=750$, (c) $t=1800$, see the text for different measures of the U-shaped domain.}\label{fig:pat1}
\end{center}
\end{figure}

Now we explain how the size of the patches and connecting channel affect the spatio-temporal pattern formation. To explain this idea, we start with the parameter values $s=3$, $a=1.65$, $d_1=0.15$ and $d_2=10$. Specific values pertaining to the domain size are $L_2=200$, $L_{x_1}=L_{x_3}=80$, $L_{x_2}=40$ and $L_y=2$ and we use the initial condition (\ref{IC1}). We present the pattern formation at different time steps in Fig. \ref{fig:pat1}. First, we observe a hot-spot pattern in the domain $\mathcal{D}_1$ and an almost uniform distribution of the populations in the rest of the domain. Then a periodic travelling wave appears in the connecting channel and spread the population in the domain $\mathcal{D}_2$ where $\mathcal{D}_2=[L_{x_1}+L_{x_2},L_1]\times[0, L_2]$. The periodic travelling wave-like pattern breaks down once it hit the boundary, and then a stationary hot spot covers the entire U-shaped domain $\mathcal{D}_2$ [see Fig.~\ref{fig:pat1}(c)]. It is interesting to note that the formation of stationary pattern in $\mathcal{D}_2$ is regulated by the width of the connecting channel (i.e. $L_y$).

\begin{figure}[ht!]
\begin{center}
        \begin{subfigure}[p]{0.31\textwidth}
                \centering
               \includegraphics[width=\textwidth]{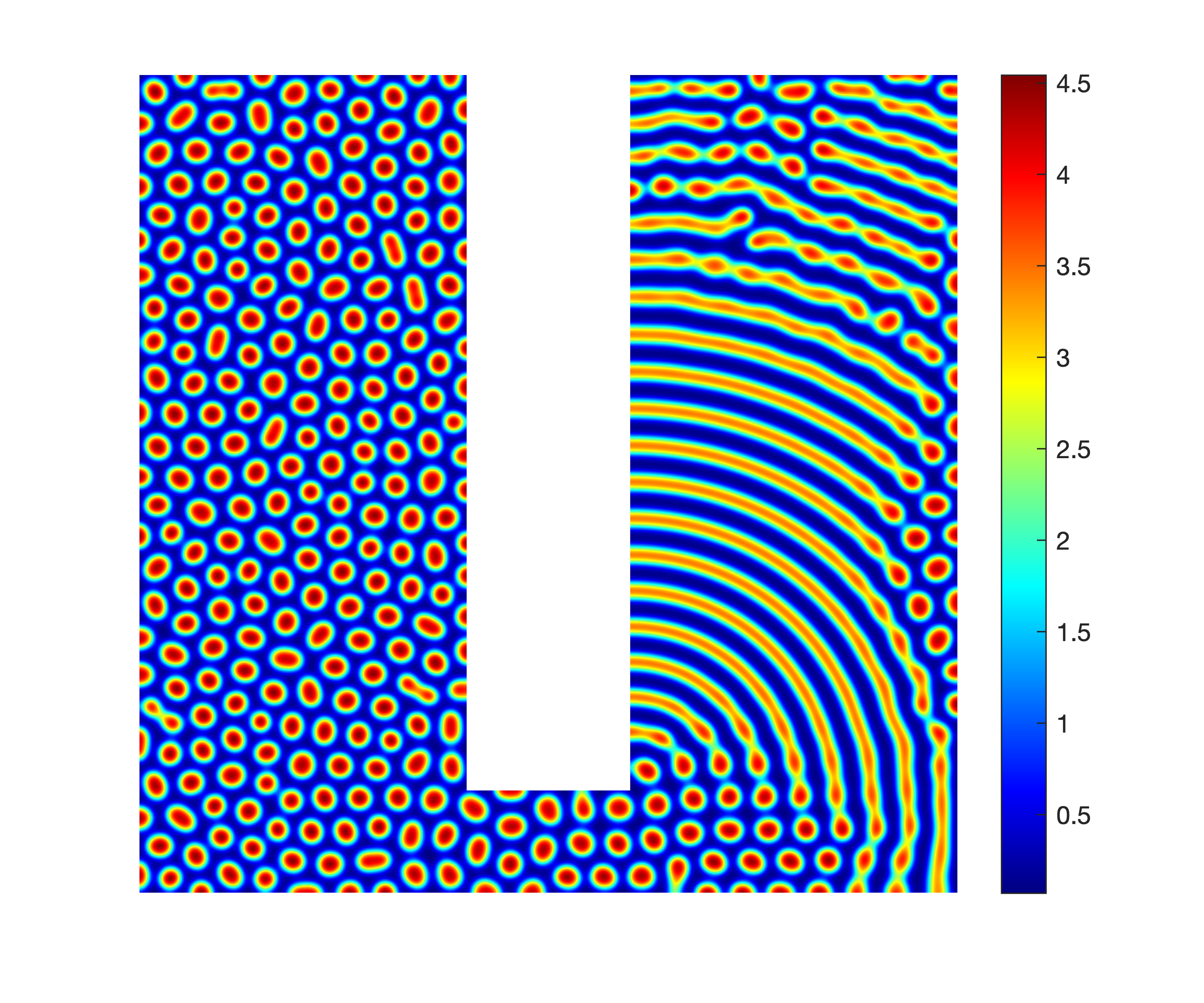}
                \caption{ }\label{fig:pat1a_transient}
        \end{subfigure}%
        ~~
        \begin{subfigure}[p]{0.31\textwidth}
                \centering
               \includegraphics[width=\textwidth]{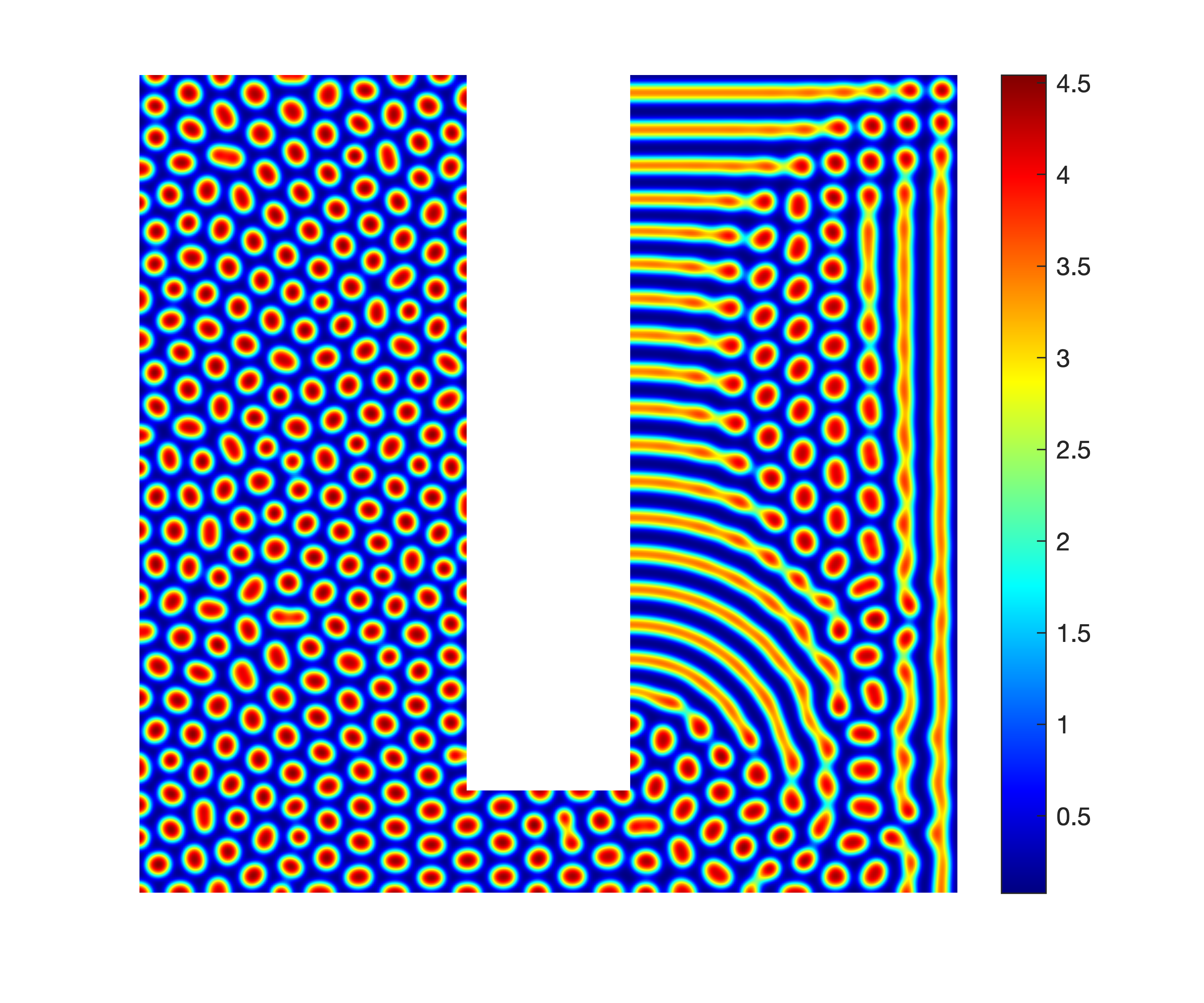}
                \caption{ }\label{fig:pat1b_transient}
        \end{subfigure}
        \caption{Transient pattern in domain $\mathcal{D}_2$ at $t=500$ for the parameter values $s=3$, $a=1.5$, $d_1=0.15$ and $d_2=10$ for two different initial conditions: (a) initial condition is (\ref{sp_IC_1}); (b) initial condition is (\ref{sp_IC_2}).}\label{fig:pat2_IC}
\end{center}
\end{figure}

\begin{figure}[htpp]
\begin{center}
        \begin{subfigure}[p]{0.31\textwidth}
                \centering
                \includegraphics[width=\textwidth]{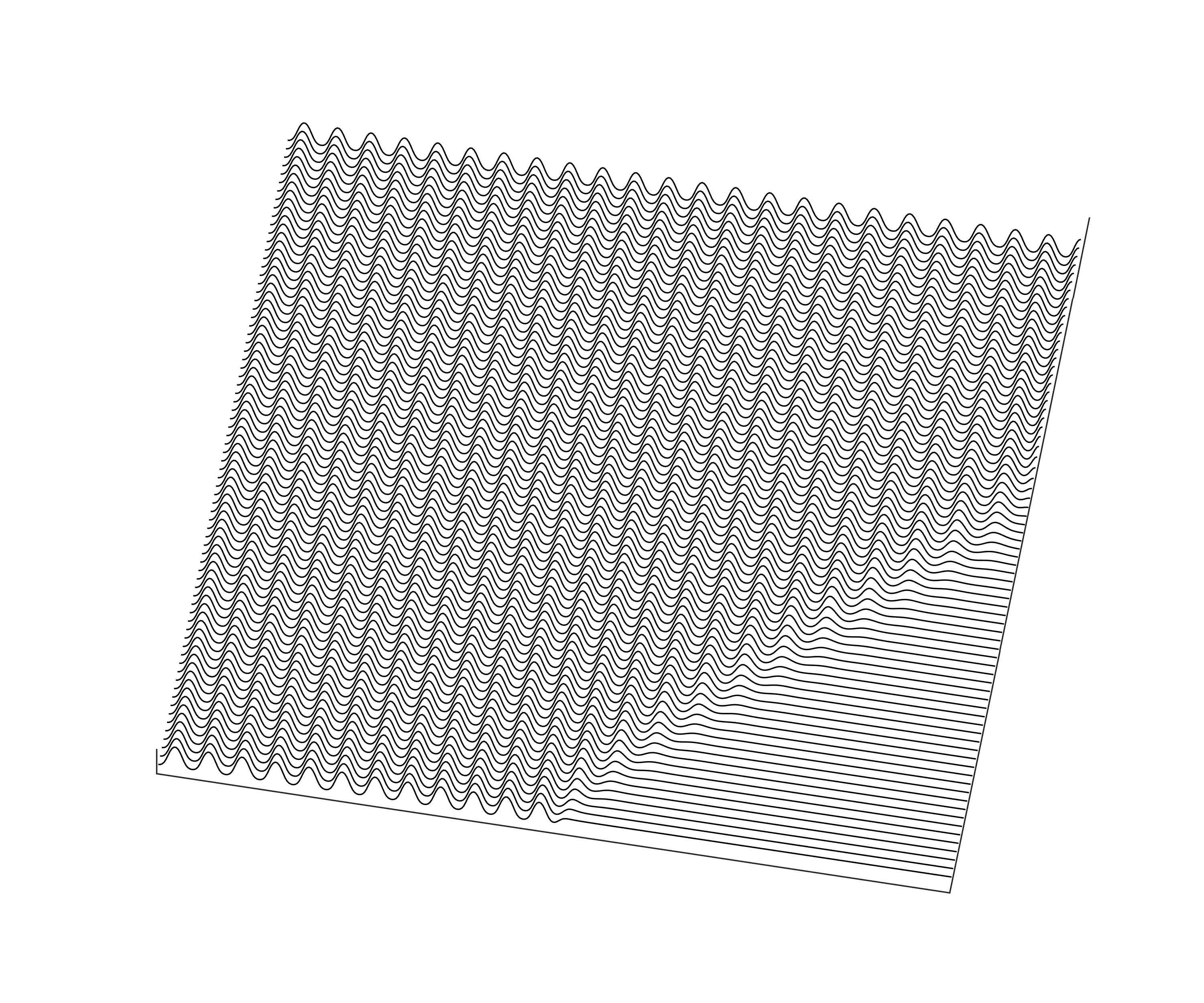}
                \caption{ }\label{fig:pat1a_invasion}
        \end{subfigure}%
        ~~~
        \begin{subfigure}[p]{0.31\textwidth}
                \centering
                \includegraphics[width=\textwidth]{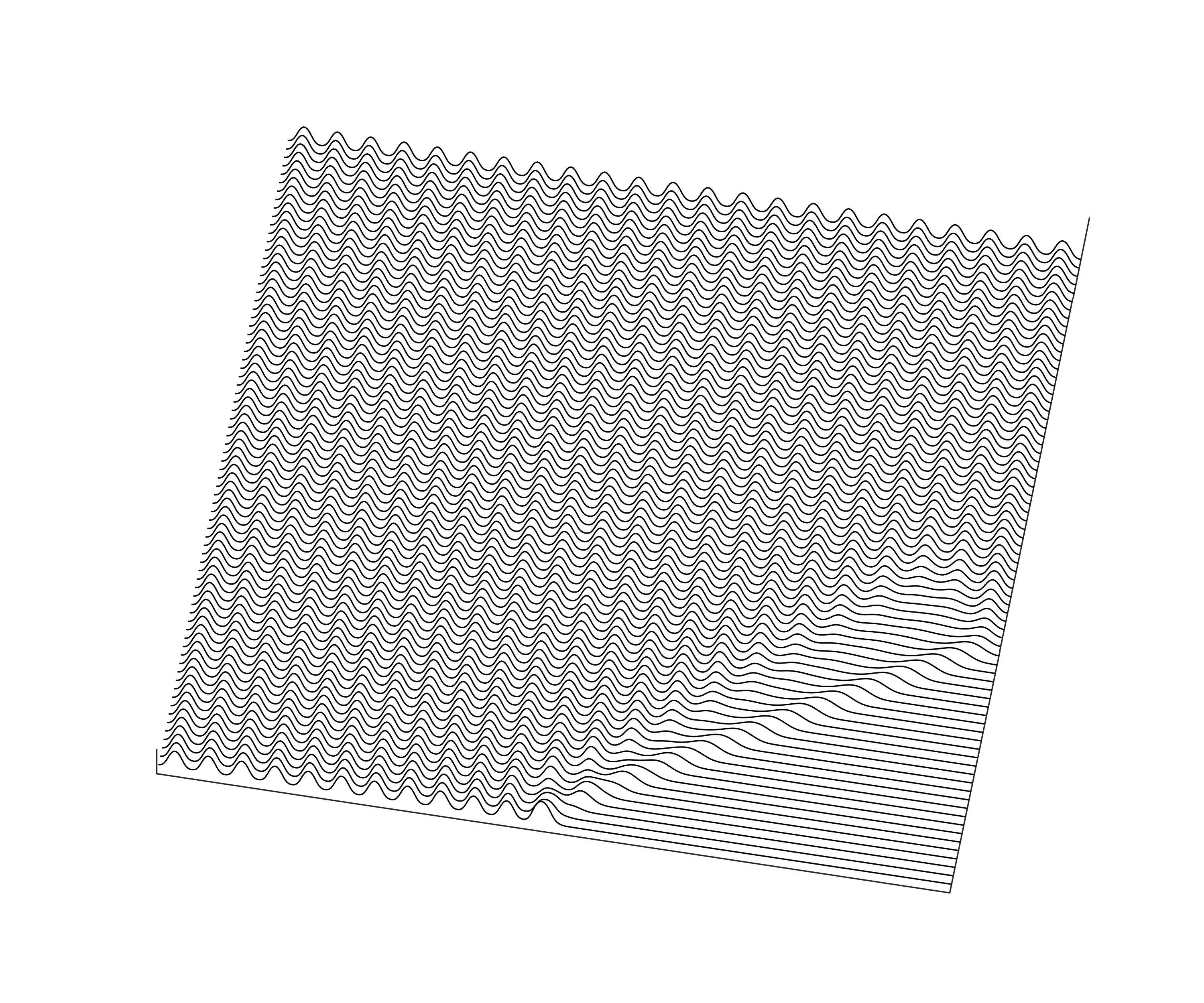}
                \caption{ }\label{fig:pat1b_invasion}
        \end{subfigure}
        \caption{Invasion of stationary hot-spot pattern formation $\mathcal{D}_1$ to $\mathcal{D}_2$ for two different initial condition in $\mathcal{D}_2$. Here $s=3$, $a=1.5$, $d_1=0.15$ and $d_2=10$, (a) $u(0,x)=u_*$, $v(0,x)=v_*$ for $x\in[L/2,L]$, (b) $u(0,x)=0$, $v(0,x)=0$ for $x\in[L/2,L]$.}\label{fig:pat:invasion}
\end{center}
\end{figure}

The initial wave propagation profile through the connecting channel to the domain $\mathcal{D}_2$ solely depend upon the width of the channel and the choice of the initial condition in $\mathcal{D}_2$. With the same choice of parameter values as mentioned in the above paragraph, we choose the initial condition in the domain $\mathcal{D}_1$ is the stationary hot-spot pattern as shown in earlier figures [cf. Fig.~\ref{fig:pat1}(c)]. We define two different initial conditions as follows:
\begin{eqnarray}\label{sp_IC_1}
u(0,x,y)=\left\{\begin{array}{ll}
\bar{u}(x,y) &  (x,y)\in\mathcal{D}_1,\\
u_* & (x,y)\in\mathcal{D}_2 \\
\end{array}\right.,
v(0,x,y)=\left\{\begin{array}{ll}
\bar{v}(x,y) &  (x,y)\in\mathcal{D}_1,\\
v_* & (x,y)\in\mathcal{D}_2 \\
\end{array}\right.,
\end{eqnarray}
and
\begin{eqnarray}\label{sp_IC_2}
u(0,x,y)=\left\{\begin{array}{ll}
\bar{u}(x,y) &  (x,y)\in\mathcal{D}_1,\\
0 & (x,y)\in\mathcal{D}_2 \\
\end{array}\right.,
v(0,x,y)=\left\{\begin{array}{ll}
\bar{v}(x,y) &  (x,y)\in\mathcal{D}_1,\\
0 & (x,y)\in\mathcal{D}_2 \\
\end{array}\right.,
\end{eqnarray}
where $\bar{u}(x,y)$ and $\bar{v}(x,y)$ is some stationary or dynamic pattern obtained through prior simulation. To be specific, we consider the stationary pattern for the prey in domain $\mathcal{D}_1$ as shown in Fig.~\ref{fig:pat1}(c) is $\bar{u}(x,y)$ and the corresponding pattern for the predator is chosen as $\bar{v}(x,y)$. We choose the U-shaped domain as $L_2=200$, $L_{x_1}=L_{x_3}=80$, $L_{x_2}=40$ and $L_y=25$. The simulation results at $t=500$ for two different initial conditions (\ref{sp_IC_1}) and (\ref{sp_IC_2}) are shown in Fig.~\ref{fig:pat2_IC} in the left and right panel, respectively. For a large time, a non-homogeneous stationary pattern occurs for both the initial conditions. But, the absence of the prey and predator species on the right side patch and the connecting channel (i.e. in $\mathcal{D}_2$) favours the rapid settlement of the stationary pattern in $\mathcal{D}_2$. We find significantly excess number of hot-spots in $\mathcal{D}_2$ corresponding to the initial condition (\ref{sp_IC_2}) while compared to the pattern obtained through the initial condition (\ref{sp_IC_1}). This is not only true for the two-dimensional case, but it is also true for the one-dimensional case. We choose a similar type of initial condition for the one-dimensional case. Here, the population distribution in the left half of the domain is the non-homogeneous stationary solution obtained through prior numerical simulation, and the right half is either a homogeneous steady-state (corresponding to the choice of parameter values) or zero population steady-state for both populations. We see that the zero population steady-state on the right half domain takes less time to settle the non-homogeneous stationary pattern over the entire domain than the homogeneous steady-state [see Fig. \ref{fig:pat:invasion}].

\begin{figure}[ht!]
        \centering
        \includegraphics[width=6cm,height=5cm]{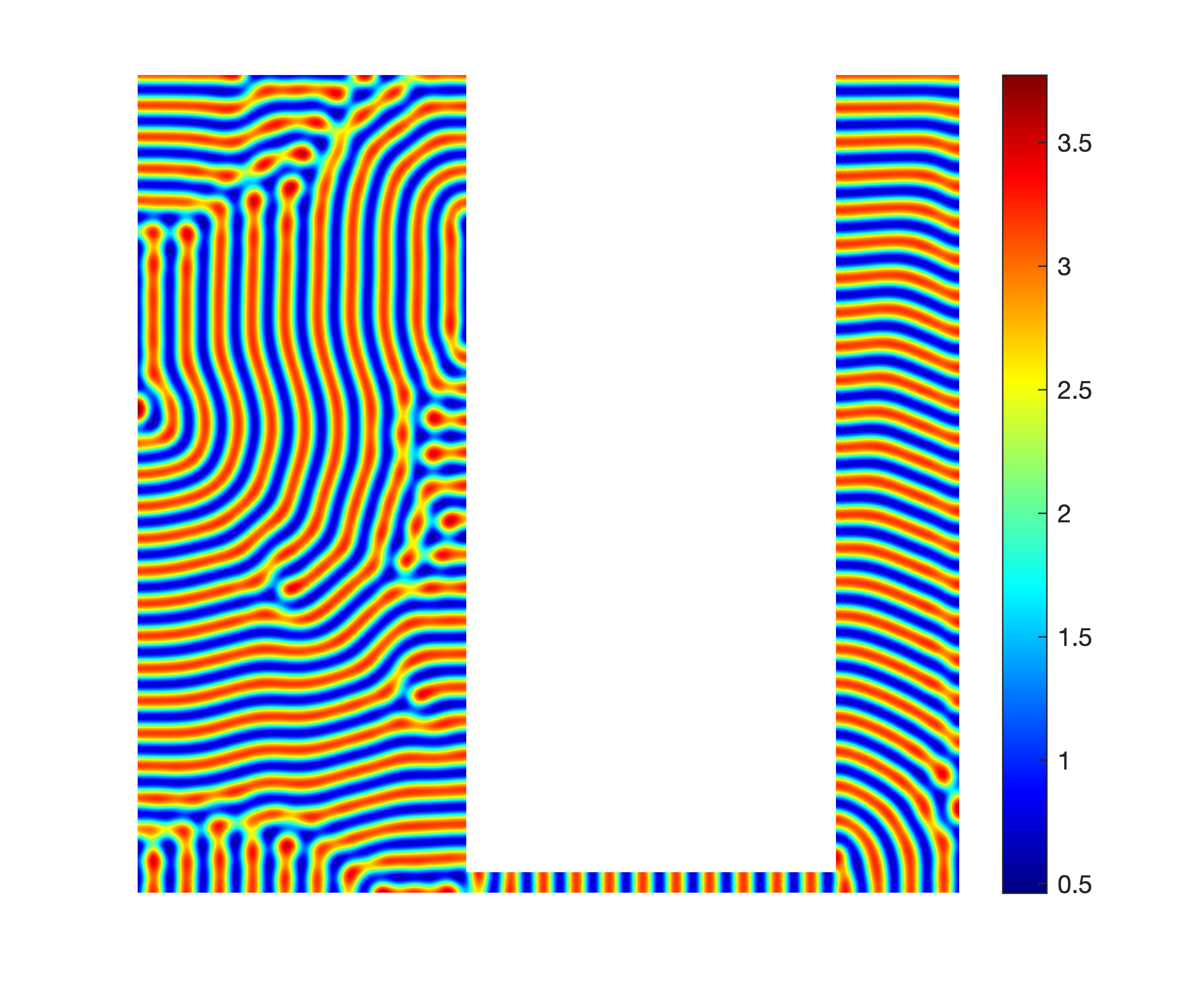}
\caption{Labyrinthine pattern in $\mathcal{D}_1$ and stripe like pattern in $\mathcal{D}_2$ for the parameters values $s=3$, $a=2$, $d_1=0.15$ and $d_2=10$.}\label{fig:pat2}
\end{figure}

The width of the connecting channel and the size of the patches are responsible for the type of stationary pattern and the time required to reach the stationary distribution. The resulting stationary pattern is related to the domain size, as the domain size allows the heterogeneous distribution to settle down to the corresponding most unstable eigenmode. To illustrate this idea, we consider the parameter values which corresponds to the labyrinthine like pattern as shown in Fig.~\ref{fig:Ushape1}(b). In order to have the right-hand patch with smaller width compared to the patch on the left, we consider the domain with $L_2=200$, $L_{x_1}=80$, $L_{x_2}=90$, $L_{x_3}=30$ and $L_y=5$. Here, we choose the initial condition as (\ref{IC1}) for the numerical simulation. We find a labyrinthine pattern in the left patch $\mathcal{D}_1$. However, a stripe pattern appears in the right patch [see Fig.~\ref{fig:pat2}]. We obtain this pattern at $t=1800$; this stationary pattern remains unaltered if we run the numerical simulations over a longer time.

\subsection{Transient pattern}

\begin{figure}[ht!]
\begin{center}
        \begin{subfigure}[p]{0.45\textwidth}
                \centering
                \includegraphics[width=\textwidth]{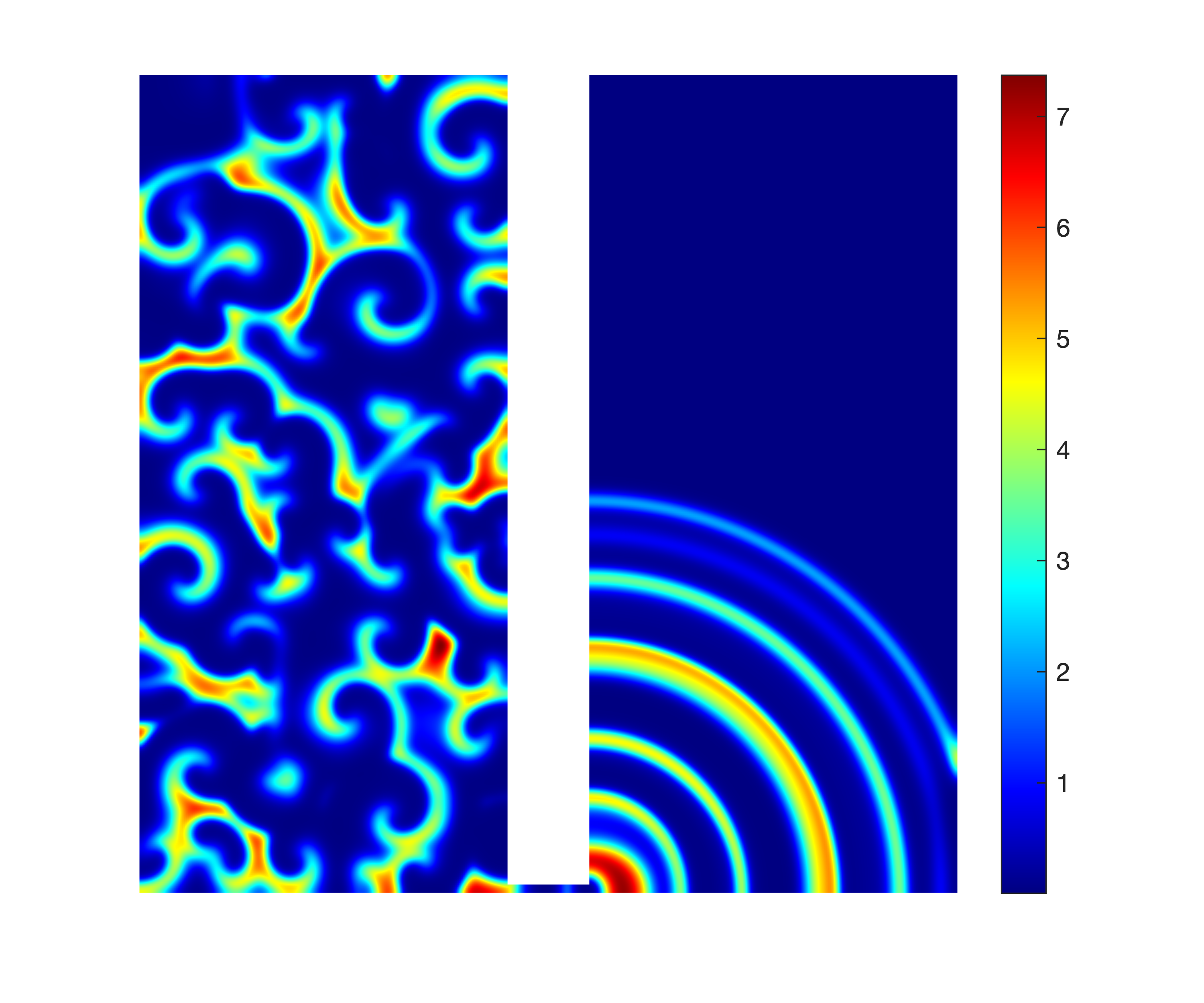}
                \caption{ }\label{fig:pat3a_transient}
        \end{subfigure}%
        ~~~
        \begin{subfigure}[p]{0.45\textwidth}
                \centering
                \includegraphics[width=\textwidth]{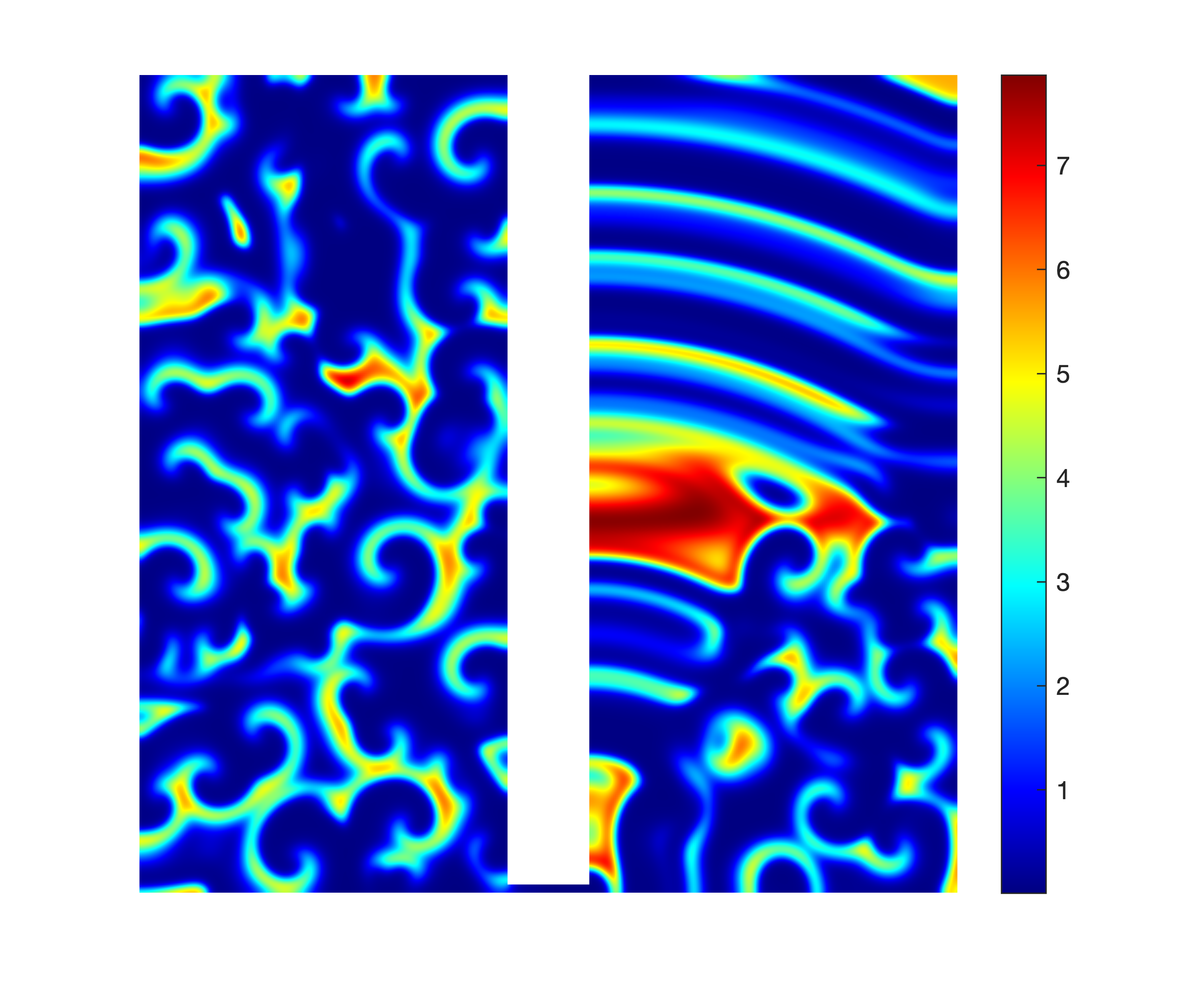}
                \caption{ }\label{fig:pat3b_transient}
        \end{subfigure}%
        \caption{Spatio-temporal chaotic pattern in $\mathcal{D}_1$ and transient pattern in $\mathcal{D}_2$. See text for parameter values, initial condition and measure of the U-shaped domain.}\label{fig:pat3}
\end{center}
\end{figure}

In order to see the formation of spatio-temporal chaotic pattern
over the U-shaped domain, we choose $a=1$, $s=2$, $d_1=d_2=1$. We
take the initial condition (\ref{IC1}) and domain size as $L_1=400$,
$L_2=200$, $L_{x_1}=L_{x_3}=180$, $L_{x_2}=40$ and $L_y=2$. We plot
the chaotic patterns at two different time steps in
Fig.~\ref{fig:pat3}. As the width of the channel is narrow, the
spatio-temporal chaotic pattern occurs on the right-hand patch after
a considerable elapsed time interval. Initially, we observe a
modulated periodic travelling wave-like pattern propagates, and then
the chaotic spatial distribution appears later. Without any loss of
generality, we can consider the modulated periodic travelling
wave-like pattern in the right-hand patch as a transient pattern. If
we look at the emerging pattern, it is challenging to realize the
periodic nature of the transient pattern. However, we can interpret
the transient pattern that appears at the right-hand patch from
initial to some advancement of time as periodic with varying
amplitude. For this, we plot the time evolution of the spatial
average of prey densities in the domains $\mathcal{D}_1$ and
$\mathcal{D}_2$.

\begin{figure}[H]
\begin{center}
\begin{subfigure}[p]{0.3\textwidth}
                \centering
                \includegraphics[width=\textwidth]{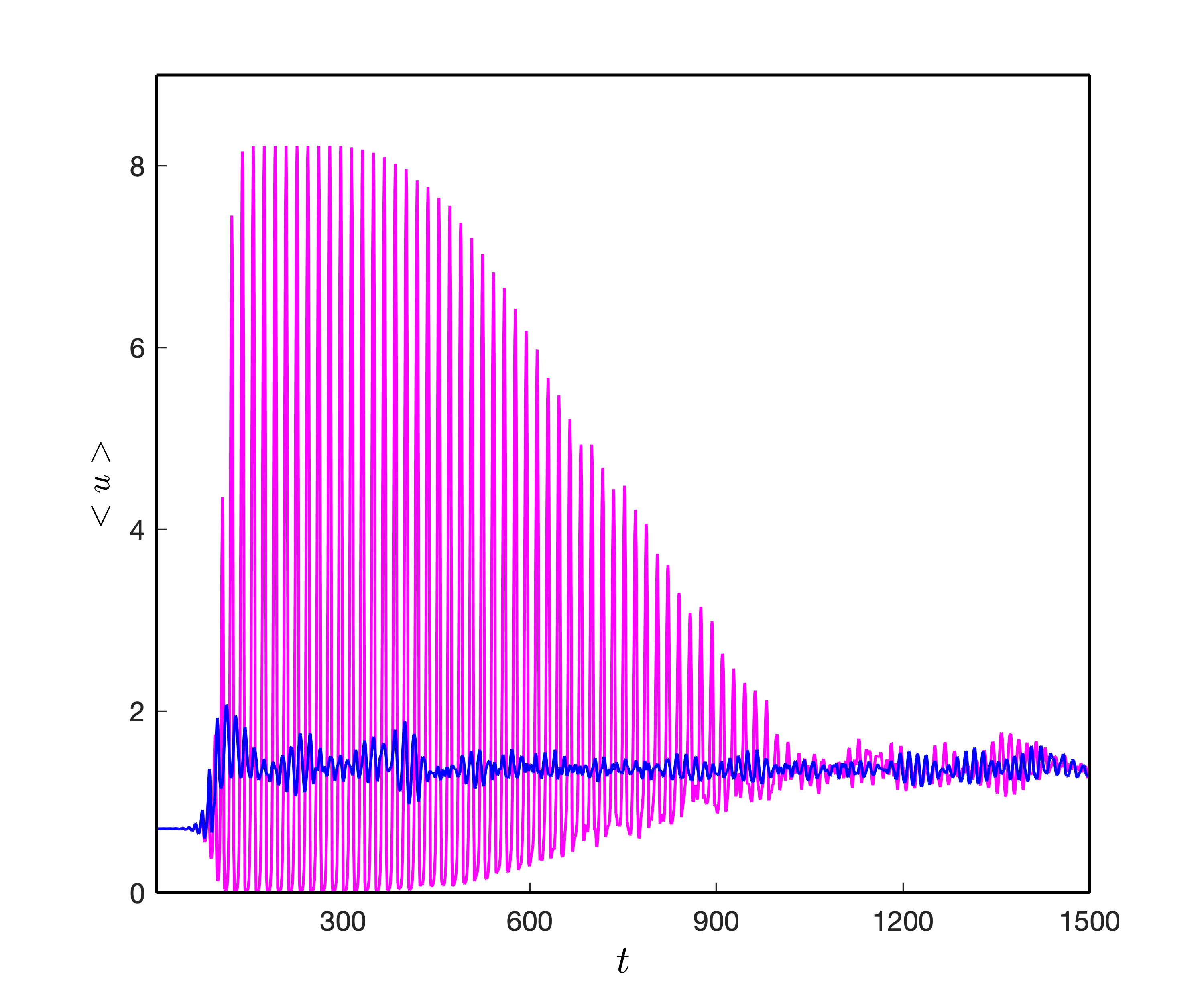}
                \caption{ }\label{fig:pat3b_L5}
        \end{subfigure}%
        \begin{subfigure}[p]{0.35\textwidth}
                \centering
                \includegraphics[width=\textwidth]{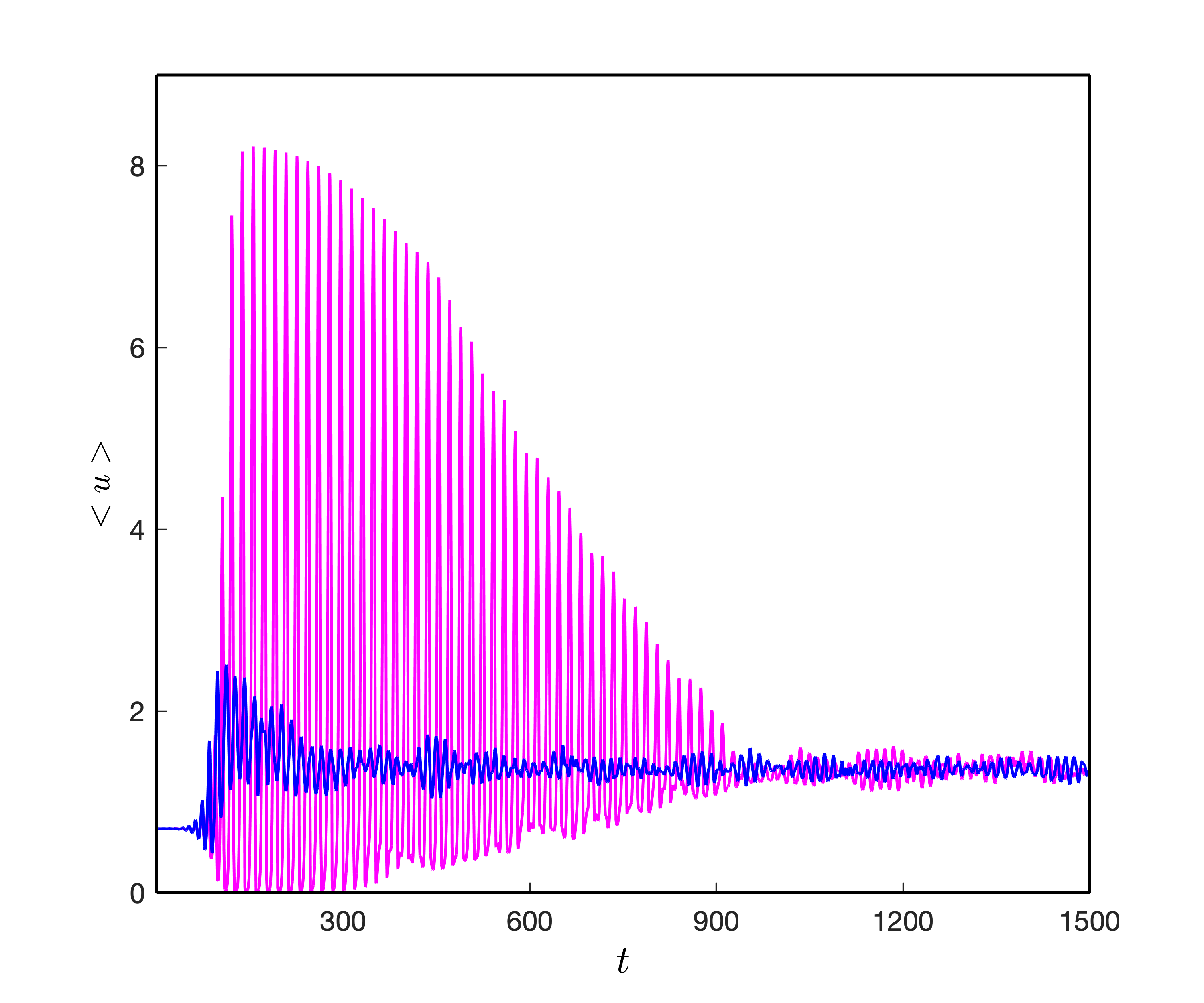}
                \caption{ }\label{fig:pat3b_L50}
        \end{subfigure}%
        \begin{subfigure}[p]{0.35\textwidth}
                \centering
                \includegraphics[width=\textwidth]{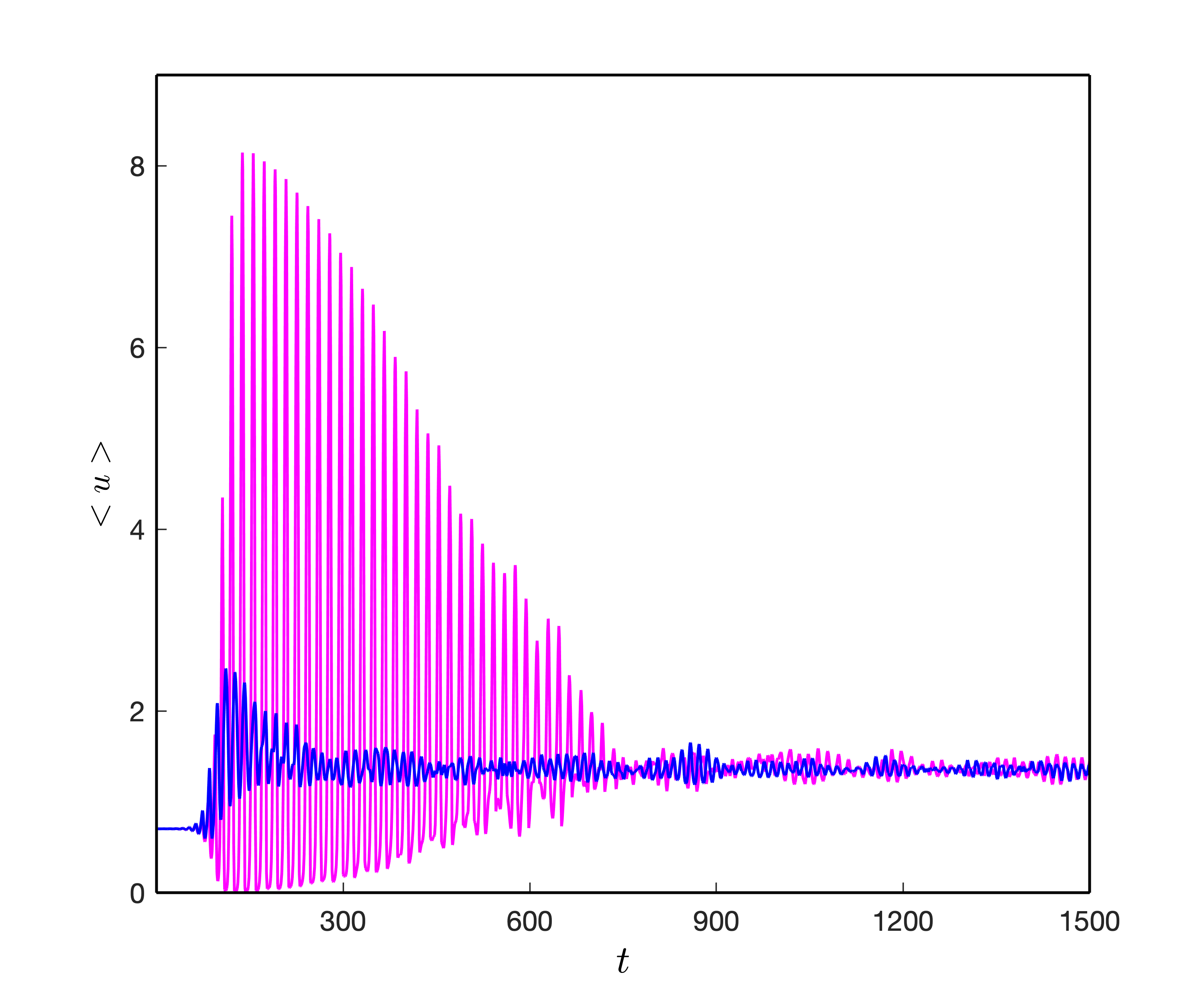}
                \caption{ }\label{fig:pat3b_L100}
        \end{subfigure}%
\end{center}
\caption{Plots of time evolution of spatial averages of prey distribution in $\mathcal{D}_1$ (blue) and $\mathcal{D}_2$ (magenta) for three different measures of the width of the connecting channel. (a) $L_y=5$, (b) $L_y=50$, (c) $L_y=100$. For other details, see the text.}\label{fig:trans1}
\end{figure}

It is interesting to report that the time evolution of spatial
average of the prey population density in $\mathcal{D}_1$ and
$\mathcal{D}_2$ exhibit significantly different behaviour during
some initial time interval. The time evolution of spatial averages
of prey density in $\mathcal{D}_1$ and $\mathcal{D}_2$ are plotted
in blue and magenta color, respectively. The numerical simulation is
performed with the initial condition (\ref{IC1}) and over a
significantly large time interval $[0,1500]$. From
Fig.~\ref{fig:trans1}, we can observe that the spatial average of
prey density in $\mathcal{D}_1$ exhibits a short transient
oscillation before entering the chaotic regime. However, the spatial
average of prey density in $\mathcal{D}_2$ shows large-amplitude
oscillation of varying magnitude for a considerable time before
entering the chaotic regime. The high amplitude periodic oscillation
of the spatial average sets in within a short time compared to the
duration over which it declines. The width of the chaotic
oscillation of the average prey density in $\mathcal{D}_2$ matches
with the average prey density in $\mathcal{D}_1$ after a
considerable time.

Finally, we like to remark that the duration of transient oscillation of average prey density in $\mathcal{D}_2$ depends on the width of the connecting channel. The time evolution of spatial averages of prey density in the domains $\mathcal{D}_1$ and $\mathcal{D}_2$ are presented in Fig.~\ref{fig:trans1} for two different values of $L_y$. We have chosen $L_y=25$ and $L_y=50$, keeping other measures related to the domain size unaltered. With an increase in the magnitude of $L_{y}$, the duration of the periodic oscillation with varying amplitudes decreases.

\section{Conclusion}

The main objective of this work is to understand the effect of fragmented habitat on spatio-temporal pattern formation. Most of the research works carried out so far, in the context of spatial pattern formation for interacting population models, is based upon the square and rectangular domain. The domain choice is mainly due to the applicability of requisite mathematical analysis and the simplicity of numerical simulations. In reality, the natural habitats of interacting populations are neither square/rectangular nor having equal width everywhere. Mostly, the boundary of the natural habitats are fractals in nature. To understand the effect of the shape and size of the fragmented habitat, we have considered the problem of spatio-temporal pattern formation in a U-shaped domain. So far as our knowledge goes, there are very few works on the pattern formation over a complex or fragmented habitat \cite{Alharbi18,Alharbi19}. Here, the U-shaped domain under consideration represents two large habitats connected by a channel through which individuals can move from one habitat to another. Earlier works \cite{Alharbi18,Alharbi19}, with an H-shaped domain, focused on the successful invasion of the exotic species from one habitat to another habitat regulated by the width of the connection between two patches. The interpretation of the fragmented domain either in the shape of `H' or in the form of `U' is the same, but the U-shaped domain is advantageous in the context of numerical simulation as it has fewer corner points. An important contribution of this work is the consideration of two types of pattern formation over the fragmented habitat, stationary and time-varying patterns.

We have considered a spatio-temporal prey-predator model with additive Allee effect in prey growth and intra-specific competition among the predators. A wide range of pattern formation scenarios by the model under consideration is already explored in detail, with both types of Allee effects - strong and weak \cite{MBWWMBE,KMMBEC,KMMBMBE}. In this work, we have concentrated within the parametric setup of the weak Allee effect only. We have explained how the strength of interaction between the species (reflected through the choice of parameter values involved with the reaction kinetics) and diffusion rate shape the stationary and dynamics patches formed by the prey and their specialist predators. As the predator species is assumed to be specialists, they always follow the prey patches to ensure the availability of their only food source. Several intrinsic factors determine the nature and size of the patches. As the growth rate of the predator is much faster than the prey, we can find several examples of such species in nature. The predator's diffusion rate compared to the prey species plays a crucial role in determining the dynamic nature of the resulting patterns. Here we have reported three different stationary patterns, namely cold-spot, labyrinthine and hot-spot, and the spatio-temporal chaotic pattern. The spatio-temporal chaotic pattern is of the form `interacting spiral pattern'. Mathematically, we have explained how the most unstable eigenmode is responsible for determining the number of stationary patches within the domain under consideration.

In the case of fragmented habitat, we have shown that the stationary and dynamic patterns can settle down over the entire domain within a reasonable time when the connecting corridor is wide enough. To explain the role of the width of the connecting corridor, we have performed numerical simulations with some specific choices of the initial condition. In most cases, we have considered a settled (non-homogeneous stationary pattern) prey and predator species in one patch. The nature of propagation of the spatial pattern to another patch depends on the width of the connecting corridor and the initial population density in the second patch. The specific choice of the initial condition described in (\ref{sp_IC_2}) helps us to understand the speed and nature of invasion of both the population in the domain $\mathcal{D}_2$. The initial population densities in the two patches and the width of the connecting corridor not only regulate the initial propagating fonts, rather has a significant impact on the transient pattern in $\mathcal{D}_2$.

Through exhaustive numerical simulations and specific choices of the width of the connecting channel, we have demonstrated various types of transient patterns in the second patch. Here we highlight two important perspectives which have not been reported so far in contemporary research works. Stationary spot patterns occur for the parameter values taken from the pure Turing domain. For such a situation, targets like periodic travelling waves can not observe. However, our numerical investigation shows that the target like patterns can emerge as a transient pattern when the width of the connecting channel is less than the width of a single stationary patch [see Fig.~\ref{fig:pat1}(a)]. For the choice of parameter values corresponding to the interacting spiral chaos, we find a wave of invasion in the second patch as a transient pattern [see Fig.~\ref{fig:pat3}(a)]. The width of the connecting channel solely regulates the persistence time of the invasive wave in the right-hand patch. Based upon these interesting observations, we can conclude that further investigation is needed to understand the spatio-temporal pattern formation over more complex/fragmented habitats. The novelty of this work is the identification of transient patterns in the context of interacting population models. Apart from the investigation of spatio-temporal pattern formation over a complex habitat, the alteration of chaotic and invasive patterns due to the presence of obstacles within the habitat is another interesting area of research \cite{JASXLCJTTNS, MJSJASNJA}. It is needless to say that research with transient phenomena in population biology is an exciting research topic nowadays \cite{AMMBSPJTB, AHscience, AHPLR, NatSci}. Identification of transient dynamics in the context of spatial pattern formation is rare in the literature; some recent works have identified quasi-steady-state transient patterns before settling to stationary state \cite{Mistro, RSI}. So far as our knowledge goes, this is the first work towards identifying transient patterns over a complex domain.

\medskip

\end{document}